\documentclass[12pt]{article}
\usepackage{bbm}
\usepackage{mathrsfs}
\usepackage{amsfonts}
\usepackage{amsmath}
\usepackage{amssymb,amsmath}
\usepackage{amscd}

\def\Z{\mathbb{Z}}
\def\C{\mathbb{C}}
\def\cl{\centerline}
\def\vs{\vspace*}

\def\QED{\hfill$\Box$}

\def\D{\Delta}
\def\A{\mathcal{A}}
\def\R{\mathcal{R}}
\def\a{\alpha}
\def\d{\cdot}
\def\c{\circ}
\def\o{\otimes}
\def\p{\partial}
\def\l{\lambda}
\def\L{\mathcal{L}(\mathcal{A})}

\def\LL{\mathcal{L}}

\numberwithin{equation}{section}
\newtheorem{theo}{Theorem}[section]
\newtheorem{defi}[theo]{Definition}
\newtheorem{exam}[theo]{Example}
\newtheorem{coro}[theo]{Corollary}

\newtheorem{prop}[theo]{Proposition}

\newtheorem{rema}[theo]{Remark}

\begin{document}

\cl{{\large\bf Super Hom-Gel'fand-Dorfman bialgebras and}}
 \cl{{\large\bf Hom-Lie
conformal superalgebras}
\footnote{Corresponding author:
lmyuan@hit.edu.cn}} \vskip8pt

\cl{ Lamei Yuan$^{\ddag}$, Sheng Chen$^{\dag}$, Caixia He$^{\dag}$} \vskip1mm
 \cl{\small{$^{\ddag}$ Academy of Fundamental and Interdisciplinary
 Sciences,}}
\cl{\small{Harbin Institute of Technology, Harbin 150080, China}}
\cl{\small{$^{\dag}$ Department of Mathematics, Harbin Institute of Technology, Harbin
150001, China}}

\cl{\small E-mail: lmyuan@hit.edu.cn, schen@hit.edu.cn, cxheedu@126.com
 }\vs{6pt}

{\small\parskip .005 truein \baselineskip 3pt \lineskip 3pt

\noindent{\bf Abstract.}  The purpose of this paper is to introduce and study super Hom-Gel'fand-Dorfman bialgebras and Hom-Lie conformal superalgebras. In this paper, we provide different ways for constructing super Hom-Gel'fand-Dorfman bialgebras and obtain some infinite-dimensional Hom-Lie superalgebras from affinization of super Hom-Gel'fand-Dorfman bialgebras. Also,
we give a general construction of Hom-Lie conformal superalgebras from Hom-Lie superalgebras and establish equivalence of quadratic Hom-Lie conformal superalgebras and super Hom-Gel'fand-Dorfman bialgebras. Finally, we characterize one-dimensional central extensions of quadratic Hom-Lie conformal superalgebras by using certain bilinear forms of super Hom-Gel'fand-Dorfman bialgebras.
 \vs{5pt}


\parskip .001 truein\baselineskip 6pt \lineskip 6pt

\vs{10pt}

\cl{\bf\S1. \
Introduction}\setcounter{section}{1}\setcounter{equation}{0}
\vs{6pt}

The Hom-type algebras arose firstly in quasi-deformation of Lie algebras of vector fields. Discrete
modifications of vector fields via twisted derivations lead to Hom-Lie algebras, which were introduced in the context of
studying deformations of Witt and Virasoro algebras \cite{HLS}.
Hom-Lie algebras were intensively discussed in a series of papers, see \cite{Yau1,Yau2,Yau3,Yau4} for example. A more general framework bordering color and super Lie algebras was introduced in \cite{LS1,LS2}. In particular, graded Hom-Lie structures including Hom-Lie super and color algebras were studied in \cite{AM,Yuan1}.

Just as Lie algebras are closely related to associative algebras, Hom-Lie algebras are
closely related to Hom-associative algebras, which were introduced in \cite{MS}. It turns out that the commutator bracket
defined by the multiplication in a Hom-associative algebra leads naturally to a Hom-Lie algebra. This provides a different way for constructing Hom-Lie algebras. Related Hom-algebras such as Hom-Novikov algebras, Hom-Poisson algebras and Hom-Lie conformal algebras have been studied in \cite{MS1,Yau,Yuan2}. Generally speaking, the notion of Hom-type algebra corresponding to some algebra $\A$ is obtained by using a linear map on $\A$ to twist the defining identities.  Yau in \cite{Yau2} gave a general method to deform an algebraic structure to the corresponding Hom-algebra via an endomorphism. This strategic way has been applied to several Hom-algebras, see \cite{MS1,Yau,Yuan2} for example.

The present paper focuses on $\Z_2$-graded Hom-algebras. Mainly, we introduce and study ``Hom-generalizations" of super Gelfand-Dorfman bialgebras and Lie conformal superalgebras, called super Hom-Gelfand-Dorfman bialgebras and Hom-Lie conformal superalgebras respectively. They are also viewed as superanalogues of Hom-Gelfand-Dorfman bialgebras and Hom-Lie conformal algebras, which we introduced and studied in \cite{Yuan2}.

In Section 2 of this paper, we summarize the main Hom-type (super)algebras and recall the notions of super Gel'fand-Dorfman bialgebra and Lie conformal superalgebra.
In Section 3, we introduce super Hom-Gel'fand-Dorfman bialgebras and provide five different ways for constructing super
Hom-Gel'fand-Dorfman bialgebras by extending some constructions of Hom-Gel'fand-Dorfman bialgebras given in \cite{Yuan2} and Hom-Novikov algebras given in \cite{Yau}. We also construct affinization of super Hom-Gel'fand-Dorfman bialgebras, which leads us to a class of infinite-dimensional Hom-Lie superalgebras. In Section 4, we introduce the notion of Hom-Lie conformal superalgebra, which is a superanalogue of Hom-Lie conformal algebras introduced in \cite{Yuan2} as well as a Hom-version of Lie conformal superalgebras. In particular, we introduce quadratic Hom-Lie conformal superalgebras and establish equivalence of quadratic Hom-Lie conformal superalgebras and super Hom-Gel'fand-Dorfman bialgebras. This generalizes the equivalent theorem in
the ungraded case obtained in \cite{Yuan2}, and in classical Lie conformal superalgebra case given in \cite{GD,X1}. Section 5 is devoted to central extensions of Hom-Lie conformal superalgebras.  We characterize one-dimensional central extensions of quadratic Hom-Lie conformal superalgebras by using certain bilinear forms of super Hom-Gel'fand-Dorfman bialgebras, generalizing a construction due to Hong \cite{H}. We should point out
that all these results could naturally be extended to more general graded cases.

Throughout this paper, all linear spaces and tensor products are over the complex field $\mathbb{C}$. In addition to the standard notation $\Z$, we use $\Z^+$ to denote the set of nonnegative
integers.

\vs{12pt}

\cl{\bf\S2. \
Preliminaries}\setcounter{section}{2}\setcounter{equation}{0}
\vs{6pt}
Let us begin with
some definitions,
including Hom-Lie superalgebra, Hom-associative algebra, Hom-Novikov superalgebra, super Gel'fand-Dorfman bialgebra, 
Hom-Gel'fand-Dorfman bialgebra, Hom-Lie conformal algebra and Lie conformal superalgebra. For
detailed discussions we refer the reader to the
literatures (e.g. \cite{AM,ka,Yau,Yuan2,ZHB} and references therein).

Let $V$ be a superspace that is a $\Z_2$-graded linear space with a direct sum
$V = V_0 \oplus V_1$. The elements of $V_j$, $j = \{0, 1\}$, are said to be homogenous and of parity $j$.
The parity of
a homogeneous element $x$ is denoted by $|x|$. Throughout what follows, if $|x|$ occurs
in an expression, then it is assumed that $x$ is homogeneous and that the expression
extends to the other elements by linearity.

\begin{defi} {\rm (see \cite{AM}) A Hom-Lie superalgebra is a superspace $\A=\A_0\oplus \A_1$ with an even
bilinear  map  $[ \cdot,\cdot]:\A\times \A\longrightarrow \A$ and an even
linear map $\alpha:\A\longrightarrow \A$, such that
\begin{eqnarray}
&&[x,y]=-(-1)^{|x||y|}[y,x],\label{super-ss}\\
&&(-1)^{|x||z|}[[x,y],\alpha(z)]+(-1)^{|x||y|}[[y,z],\alpha(x)]+(-1)^{|y||z|}[[z,x],\alpha(y)]=0,\label{super-JJ}\ \ \ \ \ \ \
\end{eqnarray}
for all homogeneous elements $x, y, z\in \A$.}
\end{defi}
\begin{rema}\rm Relation \eqref{super-ss} is called skew-symmetry, and \eqref{super-JJ} is called Hom-Jacobi identity. By \eqref{super-ss}, one may write \eqref{super-JJ} in the following form
\begin{eqnarray}\label{super-jj-j}
[\alpha(x),[y,z]]=[[x,y],\alpha(z)]+(-1)^{|x||y|}[\alpha(y),[x,z]],
\end{eqnarray}
for all homogeneous elements $x, y, z\in \A$.
\end{rema}
\begin{defi} {\rm (see \cite{AM}) A  Hom-associative superalgebra is a superspace $\A=\A_0\oplus \A_1$ with an even bilinear map
$\mu: \A \times \A \longrightarrow \A$ and an even linear map $\alpha: \A
\longrightarrow \A$, satisfying
\begin{eqnarray}\label{hom-asso}
\mu(\alpha(x),\mu(y,z)) =
\mu(\mu(x,y),\alpha(z))
\end{eqnarray}
for all homogeneous elements $x, y, z$ in $\A$.
}\end{defi}

A Hom-associative superalgebra
$(\mathcal {A},\mu,\alpha)$ is called {\it commutative }if
$\mu(x,y)=(-1)^{|x||y|}\mu(y,x)$ for $x, y\in\mathcal{A}$. A {\it derivation} on a Hom-associative superalgebra is defined in the usual way.

\begin{defi}\label{d1}{\rm  (see \cite{ZHB}) A Hom-Novikov superalgebra is a superspace $\A=\A_0\oplus \A_1$ equipped with an even bilinear
operation $\c:\A \times \A \longrightarrow \A $ and an even linear map $\alpha:\A\rightarrow \A$ such that for all homogeneous elements $x, y, z \in\A$:
\begin{eqnarray}
(x\c y)\c\alpha(z)-\alpha(x)\c(y\c z)&=&(-1)^{|x||y|}((y\c x)\c \alpha(z)
-\alpha(y)\c (x\c z)),\label{Hom-Nov1++}\\
(x\c y)\c \alpha(z)&=&(-1)^{|y||z|}(x\c z)\c \alpha(y).\label{Hom-Nov2++}
\end{eqnarray}
}\end{defi}
\begin{rema}\rm
By taking $\a={\rm id}$, we obtain Novikov superalgebras, which were studied in \cite{Xu3, Xu4}.
\end{rema}

The notion of super Gel'fand-Dorfman bialgebra was introduced and studied in \cite{X1},
where it is shown that the equivalence of super Gel'fand-Dorfman bialgebras and quadratic conformal superalgebras. This statement was essentially given in \cite{GD} (without proof).

\begin{defi}{\rm (\cite{X1}) A super Gel'fand-Dorfman bialgebra is a superspace $\A=\A_0\oplus\A_1$ equipped with two operations
$[\cdot, \cdot]$ and $\circ$, such that $(\A, [\cdot, \cdot])$ forms a Lie
superalgebra, $(\A, \circ)$ forms a Novikov superalgebra and the compatibility
condition
\begin{eqnarray}\label{GB}
[x\circ y, z]-(-1)^{|y||z|}[x\circ z, y]+[x,y]\circ z-(-1)^{|y||z|}[x,z]\circ y-x\circ[y,z]=0
\end{eqnarray}
holds for all homogeneous elements $x, y, z$ in $\A$.}
\end{defi}

The Hom-Gel'fand-Dorfman bialgebras and Hom-Lie conformal algebras were introduced and studied in \cite{Yuan2}.

\begin{defi}\label{yuan}{\rm  A Hom-Gel'fand-Dorfman
 bialgebra is a linear space $\A$ equipped with a linear selfmap $\a$ and two operations
 $[\cdot, \cdot]$ and $\circ$, such that $(\A,
[\cdot, \cdot], \a)$ is a Hom-Lie algebra, $(\A, \circ, \a)$
is a Hom-Novikov algebra and the following compatibility condition holds for $x, y, z\in\A$:
\begin{eqnarray}\label{HGB1111}
[x\circ y, \a(z)]-[x\circ z, \a(y)]+[x,y]\circ \a(z)-[x,z]\circ
\a(y)-\a(x)\circ[y,z]=0.
\end{eqnarray}}
\end{defi}

\begin{defi}\label{HLCF}{\rm
A Hom-Lie conformal algebra $\mathcal
{R}$ is a $\C[\partial ]$-module equipped with a linear selfmap $\a$ satisfying $\a\p=\p\a$, and a $\C$-bilinear map $[\cdot_\lambda\cdot]:\mathcal {R}\otimes \mathcal {R}\rightarrow
\C[\l]\otimes \mathcal{R}$ such that the following axioms
\begin{eqnarray}
[\partial a_\lambda b]&=&-\lambda[a_\lambda b],\ \ \ [a_\l \partial b]=(\partial+\l)[a_\l b],\ \  \mbox{(conformal\  sesquilinearity)}\label{HL111}\\
{[a_\lambda b]} &=& -[b_{-\lambda-\partial}a],\ \ \ \ \ \ \mbox{(skew-symmetry)}\label{HL222}\\
{[\a(a)_\lambda[b_\mu c]]}&=&[[a_\lambda b]_{\lambda+\mu
}\a(c)]+[\a(b)_\mu[a_\lambda c]]\ \ \mbox{(Hom-Jacobi \
identity)}\label{HL333}
\end{eqnarray}
hold for $a, b, c\in \mathcal {R}$.}
\end{defi}

When $\a={\rm id}$ we recover the classical Lie conformal algebra, which encodes the singular parts of the operator product expansion (OPE) of chiral fields in two-dimensional conformal field theory \cite{ka}. The $\Z_2$-graded version of Lie conformal algebras was studied in \cite{FK,ka}.

\begin{defi}\label{definition}{\rm
 A Lie conformal superalgebra $\mathcal {R}=\mathcal {R}_0\oplus \mathcal {R}_1$ is a $\C[\partial ]$-module endowed with a $\C$-bilinear map $[\cdot_\lambda\cdot]:\mathcal {R}\otimes \mathcal {R}\rightarrow
\C[\l]\otimes \mathcal{R}$,
 satisfying the following axioms ($a, b, c\in \mathcal {R}$):
 \begin{eqnarray}
[\partial a_\lambda b]&=&-\lambda[a_\lambda b],\ \ \ [a_\l \partial b]=(\partial+\l)[a_\l b], \ \ \mbox{(conformal\  sesquilinearity)}\label{L1}\\
{[a_\lambda b]} &=& -(-1)^{|a||b|}[b_{-\lambda-\partial}a],\ \ \ \ \ \ \mbox{(skew-symmetry)}\label{L2}\\
{[a_\lambda[b_\mu c]]}&=&[[a_\lambda b]_{\lambda+\mu
}c]+(-1)^{|a||b|}[b_\mu[a_\lambda c]].\ \ \mbox{(Jacobi \ identity)}\label{L3}
\end{eqnarray}
}\end{defi}
\vs{8pt}

\cl{\bf\S3. Super Hom-Gel'fand-Dorfman
bialgebras}\setcounter{section}{3}\setcounter{equation}{0}\setcounter{theo}{0}
\vs{8pt}

In this section, we introduce super Hom-Gel'fand-Dorfman bialgebras and provide five
ways for constructing super Hom-Gel'fand-Dorfman bialgebras. Also, we construct affinization of super Hom-Gel'fand-Dorfman bialgebras, leading us to a class of infinite-dimensional Hom-Lie superalgebras.

\begin{defi}{\rm  A super Hom-Gel'fand-Dorfman bialgebra (or super Hom-GD bialgebra for short) is a superspace $\A=\A_0\oplus \A_1$ equipped with an even linear map $\a:\A \rightarrow\A$ and two operations
 $[\cdot, \cdot]$ and $\circ$, such that $(\A,
[\cdot, \cdot], \a)$ is a Hom-Lie superalgebra, $(\A, \circ, \a)$
is a Hom-Novikov superalgebra and the following compatibility condition
\begin{eqnarray}\label{HGB}
[x\circ y, \a(z)]-(-1)^{|y||z|}[x\circ z, \a(y)]+[x,y]\circ \a(z)-(-1)^{|y||z|}[x,z]\circ
\a(y)-\a(x)\circ[y,z]=0
\end{eqnarray}
holds for all homogeneous elements $x, y, z \in\A$.}
\end{defi}

We recover the super Gel'fand-Dorfman bialgebras when $\a={\rm id}$.
The Hom-Gel'fand-Dorfman bialgebras are obtained when the part of parity one is trivial.

The following result shows our first construction of super Hom-Gel'fand-Dorfman bialgebras from Hom-Novikov superalgebras,  generalizing \cite[Theorem 3.2]{Yuan2}.

\begin{theo} \label{th1} Let $(\A,\circ,\a)$ be a Hom-Novikov superalgebra. Define the supercommutator on
homogeneous elements by
\begin{eqnarray}\label{commu**}
[x,y]^-=x\c y-(-1)^{|x||y|}y\c x,
\end{eqnarray}
then $(\mathcal {A},[\cdot,\cdot]^-, \circ, \alpha)$ is a super Hom-GD bialgebra.
\end{theo}
\noindent{\it Proof.~} The bracket is obviously supersymmetric. By \eqref{Hom-Nov2++} and \eqref{commu**}, we have
 \begin{eqnarray}\label{j1}
&&[x,y]^-\c \alpha(z)+(-1)^{|x|(|y|+|z|)}[y,z]^-\c \alpha(x)+(-1)^{(|x|+|y|)|z|}[z,x]^-\c \alpha(y)\nonumber\\
&=&\big((x\c y)\c \alpha(z)-(-1)^{|y||z|}(x\c z)\c \alpha(y)\big)
+(-1)^{|x||y|}\big((-1)^{|x||z|}(y\c z)\c \alpha(x)-(y\c x)\c \alpha(z)\big)\nonumber\\&&+(-1)^{(|x|+|y|)|z|}\big((z\c x)\c \alpha(y)-(-1)^{|x||y|}(z\c y)\c \alpha(x)\big)
=0.\end{eqnarray}
Similarly, we have
$$\alpha(x)\c [y,z]^-+(-1)^{|x|(|y|+|z|)}\alpha(y)\c [z,x]^-+(-1)^{|z|(|x|+|y|)}\alpha(z)\c [x,y]^-=0,$$
which together with \eqref{j1} implies \eqref{super-JJ}. Thus $(\A, [\cdot,\cdot]^-, \alpha)$ is a Hom-Lie superalgebra. By \eqref{Hom-Nov1++} and \eqref{commu**},
\begin{eqnarray*}
&&{[x \circ y, \a(z)]^- -(-1)^{|y||z|}\big([x\circ z, \a(y)]^- -[x,z]^- \circ \a(y)\big) +[x,y]^-\circ
\a(z)-\a(x)\circ [y,z]^-}\\
&=&(x \circ y)\circ \a(z)- (-1)^{(|x|+|y|)|z|}\a(z)\circ(x \circ y)-(-1)^{|y||z|}(x\circ
z)\circ\a(y)\\&&+(-1)^{|x||y|}\a(y)\circ(x\circ z)+(x\circ y)\circ\a(z)-(-1)^{|x||y|}(y\circ
x)\circ\a(z)-(-1)^{|y||z|}(x\circ z)\circ\a(y)\\&&+(-1)^{(|x|+|y|)|z|}(z\circ
x)\circ\a(y)-\a(x)\circ(y\circ z)+(-1)^{|y||z|}\a(x)\circ(z\circ y)\\
&=&\big((x\circ y)\circ\a(z) +(-1)^{|x||y|}\a(y)\circ(x\circ z)-(-1)^{|x||y|}(y\circ
x)\circ\a(z)- \a(x)\circ(y \circ z)\big)\\&&+(-1)^{|y||z|}\big((-1)^{|x||z|}(z\circ
x)\circ\a(y)+\a(x)\circ(z\circ y)- (-1)^{|x||z|}\a(z)\circ(x \circ y)-(x\circ z)\circ\a(y)\big)\\&=&0.
\end{eqnarray*}
This proves \eqref{HGB} and the theorem.\QED
\begin{exam}\label{exam2}\rm There is a three-dimensional Hom-Novikov superalgebra ($\A=\A_0\bigoplus \A_1,\circ,\alpha$) with $\A_0=\mathbb{C} x_1\bigoplus \mathbb{C} x_2$ and $\A_1=\mathbb{C} y$, such that (we only give the nonzero products)
\begin{eqnarray*}
x_1\circ x_1=x_2, \ \ x_2\circ x_2=x_1, \ \ x_1\circ y=y,\ \ x_2\circ y=y,
\end{eqnarray*}
and
\begin{eqnarray*}
 \alpha(x_1)=x_2,\ \ \! \alpha(x_2)=x_1,\ \ \ \alpha(y)=0.
\end{eqnarray*}
It is not a Novikov superalgebra since $0=(x_1\circ x_2)\circ x_1\neq(x_1\circ x_1)\circ x_2=x_1.$
By Theorem \ref{th1} we obtain a super Hom-GD bialgebra $(\A,[\cdot,\cdot]^-,\circ,\alpha)$ with the following nonzero brackets
\begin{eqnarray*}
[x_1,y]^-=y, \ [x_2,y]^-=y.
\end{eqnarray*}
\end{exam}

 The following result, given in \cite{ZHB}, extends Yau' construction of Hom-Novikov algebras from Novikov algebras with an algebra endomorphism \cite{Yau}.
\begin{prop} \label{prop2+2} Let $(\mathcal
{A}, \circ)$ be a Novikov superalgebra with an even algebra homomorphism $\a$. Then $(\A,\circ_\alpha,\alpha)$ is a
Hom-Novikov superalgebra with
\begin{eqnarray}
x \circ_{\a} y =\alpha(x)\circ\alpha(y),\  \forall\ x, y\in\A.
\end{eqnarray}
\end{prop}

 As a consequence of Theorem \ref{th1} and Proposition
\ref{prop2+2}, we have
\begin{coro}\label{co} Let $(\mathcal
{A}, \circ)$  be a Novikov superalgebra with an even algebra endomorphism $\a$. Then $(\mathcal {A}, [\cdot,\cdot]^-_{\a},\circ_{\a}, \a)$ is a super Hom-GD bialgebra, such that
\begin{eqnarray}
[x,y]^-_{\a}=\a(x\circ y)-(-1)^{|x||y|}\a(y\circ x),
\end{eqnarray}
for all homogeneous elements $x, y$ in $\A$.
\end{coro}

\begin{exam}\label{exam3}\rm Consider the following three-dimensional Novikov superalgebra $\A=\A_0\mbox{$\bigoplus$} \A_1$ with $\A_0=\mathbb{C}x_1\bigoplus \mathbb{C}x_2$, $\A_1=\mathbb{C}y$, satisfying the following nontrivial products
\begin{eqnarray}
 x_1\circ x_2=\frac12x_1,\ x_2\circ x_1=-\frac12x_1,\ x_2\circ x_2=\frac12x_2,\ y\circ x_2=\frac12y,\ y\circ y=\frac12x_1.
\end{eqnarray}
For any fixed $\lambda \in \mathbb{C}$, define an algebra homomorphism $\alpha : \A\longrightarrow \A$ by
\begin{eqnarray*}
\alpha (x_1)=\lambda^2x_1,\ \alpha (x_2)=x_2,\ \alpha (y)=\lambda y.
\end{eqnarray*}
By Proposition
\ref{prop2+2} we get a Hom-Novikov superalgebra $(\A,\circ_\alpha,\alpha)$
with the following nonzero products
\begin{eqnarray}
x_1\circ_\alpha x_2=\frac{\lambda^2}{2} x_1,\ x_2\circ_\alpha x_1=-\frac{\lambda^2}{2}x_1,\ x_2\circ_\alpha x_2=\frac12x_2,\ y\circ_\alpha x_2=\frac{\lambda}{2}y,\ y\circ_\alpha y=\frac{\lambda^2}{2} x_1.
\end{eqnarray}
Note that it is not a Novikov superalgerbra when $\lambda\neq 0,1$.
According to Corollary \ref{co} we obtain a super Hom-GD bialgebra $(\A,[\cdot,\cdot]^-_{\a},\circ_\a,\alpha)$ with the following nonzero products
\begin{eqnarray*}
[x_1,x_2]^-_{\a}=\lambda^2x_1, \ [x_2,y]^-_{\a}=-\frac\lambda 2 y.
\end{eqnarray*}
\end{exam}

Our second construction of super Hom-GD bialgebras comes from
super Gel'fand-Dorfman bialgebras along with an even algebra endomorphism. It is a superanalogue of \cite[Theorem 3.5]{Yuan2}.
\begin{theo}\label{theo-2}
Let $(\mathcal {A}, [\cdot,\cdot],\circ)$ be a super Gel'fand-Dorfman
algebra equipped with an even algebra endomorphism $\alpha$. Then $(\mathcal {A},
[\cdot,\cdot]_\a,\circ_\a, \a)$ is a super Hom-GD bialgebra with
\begin{eqnarray}
x \circ_{\a} y =\alpha(x)\circ\alpha(y), \ \ \ [x,y]_\a=[\a(x),\a(y)], \ \forall\ x, y\in\A.
\end{eqnarray}
\end{theo}
\noindent{\it Proof.~} It is known that $(\mathcal {A}, [\cdot,\cdot]_\a, \a)$
is a Hom-Lie superalgebra and $(\mathcal {A}, \circ_\a, \a)$ is a
Hom-Novikov superalgebra. By the fact that $\a$ is an even algebra endomorphism of $(\mathcal {A}, [\cdot,\cdot],\circ)$ and \eqref{GB},
\begin{eqnarray*}
&&[x
\circ_\a y, \a(z)]_\a -(-1)^{|y||z|}\big([x\circ_\a z, \a(y)]_\a-[x,z]_\a \circ_\a \a(y)\big) +[x,y]_\a\circ_\a
\a(z) - \a(x)\circ_\a [y,z]_\a\\
&=& \a^2\big([x \circ y, z] -(-1)^{|y||z|}[x\circ z, y] +[x,y]\circ
z-(-1)^{|y||z|}[x,z] \circ y- x\circ [y,z]\big)=0.
\end{eqnarray*} This concludes the proof.
 \QED\vskip8pt

\begin{exam}\rm Let $(\A,\cdot)$ be a commutative associative algebra with a derivation $d$ and an algebra homomorphism $\alpha$, such that $\alpha d=d \alpha$. Set
\begin{eqnarray*}
\widetilde{\A}=\A \times \A =\widetilde{\A_0} \mbox{$\bigoplus$} \widetilde{\A_1},\ \mbox{with}\ \widetilde{\A_0}=(\A,0),  \ \widetilde{\A_1}=(0,\A).
\end{eqnarray*}
For any $u_0,  u_1,  v_0,  v_1 \in \A$, define two algebraic operations $[\cdot,\cdot]$ and $\circ$ on $\widetilde{\A}$ by
\begin{eqnarray}
[(u_0,u_1),(v_0,v_1)]&=&(u_1\cdot v_1,0),\\
(u_0,u_1)\circ (v_0,v_1)&=&(u_0\cdot d(v_0),u_1\cdot d(v_0)+u_0\cdot d(v_1)).
\end{eqnarray}
They define a super Gel'fand-Dorfman bialgebra $(\widetilde{\A},[\cdot,\cdot],\circ)$ by \cite[Theorem 3.3]{X1}.

Define a linear map $\tilde \a:\widetilde{\A}\rightarrow \widetilde{\A}$ by
\begin{eqnarray}
\tilde\alpha (u,v)=(\alpha (u),\alpha (v)), \ \forall \ u,v\in\A.
\end{eqnarray}
For any $u_0,  u_1,  v_0,  v_1 \in \A$, we have
\begin{eqnarray*}
\tilde\alpha [(u_0,u_1),(v_0,v_1)]&=&\tilde\alpha (u_1\cdot v_1,0)=(\alpha (u_1\cdot v_1),0)=(\alpha (u_1)\cdot \alpha (v_1),0)\nonumber\\
&=&[(\alpha (u_0),\alpha (u_1)),(\alpha (v_0),\alpha (v_1))]\nonumber\\
&=&[\tilde\alpha (u_0,u_1),\tilde\alpha (v_0,v_1)],\\
\tilde\alpha ((u_0,u_1)\circ (v_0,v_1))&=&\tilde\alpha (u_0\cdot d(v_0),u_1\cdot d(v_0)+u_0\cdot d(v_1))\\
&=&(\alpha (u_0)\cdot \alpha d(v_0),\alpha (u_1)\cdot \alpha d(v_0)+\alpha (u_0)\cdot \alpha d(v_1))\\
&=&(\alpha (u_0)\cdot d\alpha (v_0),\alpha (u_1)\cdot d\alpha(v_0)+\alpha (u_0)\cdot d\alpha (v_1))\\
&=&(\alpha (u_0),\alpha (u_1))\circ (\alpha (v_0),\alpha (v_1))\\
&=&\tilde\alpha (u_0,u_1)\circ \tilde\alpha (v_0,v_1),
\end{eqnarray*}
since $\alpha d=d\alpha$. Thus $\tilde\alpha$ is an even algebra endomorphism.
 By Theorem \ref{theo-2}, we obtain a super Hom-GD bialgebra $(\widetilde{\A},[\cdot,\cdot]_{\tilde\a},\circ_{\tilde\a},\tilde\a)$
with
\begin{eqnarray*}
[(u_0,u_1),(v_0,v_1)]_{\tilde\a}&=&(\a(u_1\cdot v_1),0),\\
(u_0,u_1)\circ (v_0,v_1)_{\tilde\a}&=&(\alpha (u_0)\cdot d\alpha (v_0),\alpha (u_1)\cdot d\alpha(v_0)+\alpha (u_0)\cdot d\alpha (v_1)),
\end{eqnarray*}
for all $u_0,  u_1,  v_0,  v_1 \in \A$.
\end{exam}

 Let $(\A,\cdot)$ be a commutative associative algebra endowed with a derivation $D$. It is proved in \cite{F,GD,Xu1} that the following operation ($\lambda$ is a fixed element in $\C$ or $\A$)
 \begin{eqnarray}\label{Gd1}
x\circ y=x\d D(y)+\lambda x\d y, \ \forall\
x, y\in\mathcal {A},
\end{eqnarray}
equips $\mathcal {A}$ with a Novikov algebra structure. Yau in \cite{Yau} gave an analogous construction of Hom-Novikov algebras in the case when $\lambda=0$. In the following we extend this construction to Hom-Novikov superalgebras.

\begin{prop}\label{twist111}
Let $(\mathcal {A},\d, \alpha)$ be a commutative Hom-associative
superalgebra with an even derivation $D$ such that $\alpha D=D\a$. Let $\lambda$
be a fixed complex number. Define
\begin{eqnarray}\label{Gd1111}
x\circ y=x\d D(y)+\lambda x\d y, \ \forall\
x, y\in\mathcal {A},
\end{eqnarray}
Then
$(\mathcal {A},\circ, \alpha)$ is a Hom-Novikov superalgebra.
\end{prop}
\noindent{\it Proof.~} By \eqref{hom-asso}, \eqref{Gd1111} and the fact that $\alpha D=D\a$, we have
\begin{eqnarray}
(x\c y)\c\a(z)&=&(x\d D(y)+\lambda x\d y)\c\a(z)\nonumber\\
&=&(x\d D(y)+\lambda x\d y)\d D(\a(z))+\lambda (x\d D(y)+\lambda x\d
y)\d \a(z)\nonumber\\&=& \a(x)\d (D(y)\d D(z))+\lambda \a(x)\d
\big(y\d D(z)+D(y)\d z\big)+\lambda^2 \a(x)\d (y\d z). \ \ \ \ \label{1}
\end{eqnarray}
 It immediately follows that
\begin{eqnarray*}
(x\c z)\c\a(y)=\a(x)\d (D(z)\d D(y))+\lambda\a(x)\d\big(z\d
D(y)+D(z)\d y\big)+\lambda^2 \a(x)\d (z\d y),
\end{eqnarray*}
which combines with \eqref{1} gives $$(x\c y)\c\a(z)=(-1)^{|y||z|}(x\c z)\c\a(y),$$
since $(\A,\d)$ is commutative and $D$ is an even derivation. Thus \eqref{Hom-Nov2++} holds. By commutativity and Hom-associativity of $(\A,\d,\a)$,
\begin{eqnarray*}
(x\c y)\c\a(z)-\a(x)\c(y\c z)
&=&-(x\d y) \d \a(D^2(z))-\lambda(x\d y) \d \a(D(z))\\
&=&-(-1)^{|x||y|}\big((y\d x) \d \a(D^2(z))+\lambda(y\d x) \d \a(D(z))\big)\\
&=&(-1)^{|x||y|}\big((y\c x)\c\a(z)-\a(y)\c(x\c z)\big),
\end{eqnarray*}
which proves \eqref{Hom-Nov1++}. Thus $(\mathcal {A}, \circ, \alpha)$
is a Hom-Novikov superalgebra.\QED\vskip8pt

The following result, which is a $\Z_2$-graded version of \cite[Theorem 3.8]{Yuan2},
shows our third construction of super Hom-GD bialgebras
from a commutative Hom-associative superalgebra with an even derivation.
\begin{theo}\label{twist}
Let $(\mathcal {A},\d, \alpha)$ be a commutative Hom-associative
superalgebra with an even derivation $D$ such that $\alpha D=D\a$. Then
$(\mathcal {A},[\d,\d]^-, \circ, \alpha)$ is a super Hom-GD bialgebra,
where $\circ$ is defined by (\ref{Gd1111}) and
\begin{eqnarray}
[x,y]^-=x\d D(y)-(-1)^{|x||y|}y\d D(x), \ \forall \ x,y\in\mathcal{A}.
\end{eqnarray}
\end{theo}
\noindent{\it Proof.~} It follows from Proposition \ref{twist111} and Theorem \ref{th1}.\QED\vskip8pt

\begin{exam}\label{exam33333}\rm Let $(\A,\cdot,\alpha)$ be a commutative Hom-associative algebra with a derivation $d$ such that $\alpha d=d \alpha$. Set
\begin{eqnarray*}
\widetilde{\A}=\A \times \A =\widetilde{\A_0} \mbox{$\bigoplus$} \widetilde{\A_1}, \ \mbox{with} \ \widetilde{\A_0}=(\A,0),  \ \widetilde{\A_1}=(0,\A).
\end{eqnarray*}
For any $u_0, u_1,  v_0,  v_1 \in \A$, define an algebraic operation $\star$ on $\widetilde{\A}$ by
\begin{eqnarray}
(u_0,u_1)\star (v_0,v_1)=(u_0\cdot v_0,0),
\end{eqnarray}
and two even linear maps $\a:\widetilde{\A}\rightarrow \widetilde{\A}$ and $D:\widetilde{\A}\rightarrow \widetilde{\A}$ by
\begin{eqnarray}
\tilde\alpha (u_0,u_1)=(\alpha (u_0),\alpha (u_1)),\ \ D(u_0,u_1)=(d(u_0),u_1).
\end{eqnarray}
 Clearly, $(\widetilde{\A}, \star,\tilde\a)$ is a commutative Hom-associative superalgebra, and $\tilde\a D= D\tilde\a$
 since $\alpha d=d\alpha$. Because $d$ is a derivation of $(\A,\cdot)$, we have
\begin{eqnarray}\label{dd1}
D((u_0,u_1)\star (v_0,v_1))=D (u_0\cdot v_0,0)=(d(u_o\cdot v_0),0)=(d(u_0)\cdot v_0+u_0\cdot d(v_0),0),
\end{eqnarray}
and
\begin{eqnarray}
D(u_0,u_1)\star (v_0,v_1)+(u_0,u_1)\star D(v_0,v_1)&=&(d(u_0),u_1)\star (v_0,v_1)+(u_0,u_1)\star (d(v_0),v_1)\nonumber\\
&=&(d(u_0)\cdot v_0,0)+(u_0\cdot d(v_0),0)\nonumber\\
&=&(d(u_0)\cdot v_0+u_0\cdot d(v_0),0).\label{dd2}
\end{eqnarray}
Combining \eqref{dd1} with \eqref{dd2} gives
\begin{eqnarray}
D((u_0,u_1)\star (v_0,v_1))=D(u_0,u_1)\star (v_0,v_1)+(u_0,u_1)\star D(v_0,v_1).
\end{eqnarray}
 Thus $D$ is a derivation of $(\widetilde{\A},\star)$. For any fixed complex number $\lambda$, the new operations
\begin{eqnarray}
x\circ y=x\star D(y)+\lambda x\star y,\ \
{[x,y]^-}=x\star D(y)-(-1)^{|x||y|}y\star D(x)
\end{eqnarray}
define a super Hom-GD bialgebra $(\tilde\A,[\cdot,\cdot]^-,\circ,\tilde\a)$ by Theorem \ref{twist}.
\end{exam}
\begin{coro}\label{coro} Let $(\A,\d)$ be a commutative associative superalgebra equipped with
an even algebra endomorphism $\a$ and
an even derivation $D$ satisfying $D\a=\a D$. For any fixed complex number $\lambda$, define
\begin{eqnarray}\label{gd3}
x\circ y=\a(x\d D(y))+\lambda\a(x\d y), \ \forall\
x,y\in\mathcal {A}.
\end{eqnarray}
Then $(\A,[\d,\d]^-, \c,\a)$ is a super Hom-GD bialgebra, where
\begin{eqnarray}
[x,y]^-=\a(x\d D(y))-(-1)^{|x||y|}\a(y\d D(x)), \ \forall \ x,y\in\mathcal {A}.
\end{eqnarray}
\end{coro}
\noindent{\it Proof.~} Write $x\d_\a y=\a(x\d y)$ for
$x,\ y\in\A$. Because $\a$ is an algebra endomorphism,
$(\A,\d_\a,\a)$ is a commutative Hom-associative superalgebra, and $D$ is also a derivation of $(\A,\d_\a)$ since
$\a D=D\a$. Rewrite
\eqref{gd3} as
\begin{eqnarray*}
x\circ y=x\d_\a D(y)+\lambda x\d_\a y,\ \forall\
x,y\in\mathcal {A}.
\end{eqnarray*}
According to Theorem \ref{twist}, $(\A,[\d,\d]^-, \c,\a)$ is a
 super Hom-GD bialgebra.\QED\vskip8pt

Our fourth construction of super Hom-GD bialgebras is related to
the following concept.

\begin{defi}{\rm A Hom-Poisson superalgebra is  a superspace $\A=\A_0\oplus \A_1$ equipped with
two operations $\cdot$ and $[\cdot,\cdot]$, and a linear
 selfmap $\a$, such that $(\A,\cdot,\a)$ is a commutative Hom-associative superalgebra,
$(\A,[\cdot,\cdot],\a)$ is a Hom-Lie superalgebra, and the following relation
\begin{eqnarray}\label{p}
[\a(x),y\d z] = (-1)^{|x||y|}\a(y)\d[x, z] + (-1)^{|z|(|x|+|y|)}\a(z)\d[x, y]
\end{eqnarray}
holds for all homogeneous elements $x, y, z\in \A$.}
\end{defi}

By setting $\A=\{0\}$ we recover the notion of
Hom-Poisson algebra, which was introduced
in the study of deformation theory of Hom-Lie algebras \cite{MS1}.
Relation \eqref{p} can be reformulated as
\begin{eqnarray}\label{p-1}
[x\d y,\a(z)] = (-1)^{|y||z|}[x, z]\d\a(y)+ \a(x)\d[y, z],
\end{eqnarray}
for all homogeneous elements $x, y, z\in \A$.
\begin{theo} \label{th2} Let $(\A, \d, [\cdot,\cdot],\a)$ be a Hom-Poisson
superalgebra. Let $D$ be an even derivation of $(\A, \d)$ such that $D\a=\a D$,
and
\begin{eqnarray}\label{gd2}
D([x,y])=[D(x),y]+[x,D(y)]+\lambda[x,y],\ \forall \ x,y\in\A,
\end{eqnarray}
where $\lambda$ is a fixed complex number. Define
\begin{eqnarray}\label{hhh} x\circ y=x\d
D(y)+\lambda x\d y, \ \forall \ x,y\in\A.
\end{eqnarray}
Then $(\mathcal {A},[\d,\d],\circ,\a)$ is a super Hom-GD bialgebra.
\end{theo}
\noindent{\it Proof.~} By Proposition \ref{twist111}, $(\mathcal
{A},\circ,\a)$ is a Hom-Novikov superalgebra. By \eqref{p-1}--\eqref{hhh} and the fact that $D$ is an even derivation of
$(\A,\d)$ commuting with $\a$, we have
\begin{eqnarray*}
&&{[x \circ y, \a(z)] -(-1)^{|y||z|}[x\circ z, \a(y)] +[x,y]\circ
\a(z)-(-1)^{|y||z|}[x,z] \circ \a(y)- \a(x)\circ [y,z]}\\
& =& [x\d D(y)+\lambda x\d y, \a(z)]-(-1)^{|y||z|}[x\d D(z)+\lambda x\d
z,\a(y)]+[x,y]\d \a(D(z))+\lambda [x,y]\d\a(z)\\&& -(-1)^{|y||z|}[x,z]\d
\a(D(y))-\lambda(-1)^{|y||z|} [x,z]\d\a(y)-\a(x)\d D([y,z])-\lambda\a(x)\d [y,z]
\\&=&[x\d D(y), \a(z)]-(-1)^{|y||z|}[x\d D(z)+\lambda x\d
z,\a(y)]+[x,y]\d \a(D(z))+\lambda [x,y]\d\a(z)\\&& -(-1)^{|y||z|}[x,z]\d
\a(D(y))-\a(x)\d \big([D(y),z]+[y,D(z)]+\lambda[y,z]\big)\\
&=&\big([x\d
D(y),\a(z)]-(-1)^{|y||z|}[x,z]\d\a(D(y))-\a(x)\d[D(y),z]\big)\\&& -\big((-1)^{|y||z|}[x\d
D(z),\a(y)]-[x,y]\d\a(D(z))-(-1)^{|y||z|}\a(x)\d[D(z),y]\big)\\&& +\lambda\big(-(-1)^{|y||z|}[x\d z,\a(y)]+[x,y]\d\a(z)-\a(x)\d[y,z]\big)=0,
\end{eqnarray*}
which proves \eqref{HGB} and the theorem. \QED\vskip8pt

\begin{exam}\rm Let $(\A, \cdot, [\cdot, \cdot], \alpha)$ be a three-dimensional Hom-Poisson superalgebra, where
$\A_{0}=\C e_1\oplus\C e_2 $ and $\A_{1}=\C e_3$, and $\a$ is an algebra endomorphism
defined by 
$$\alpha(e_{1})=ae_{1},\ \
\alpha(e_{2})=e_{1}+e_{2}, \ \ \alpha(e_{3})=\mu e_{3},$$
where $a$ and $\mu$ are fixed nonzero complex numbers.
The defining relations (we give only the ones with nonzero values in the right-hand side) are
$$e_{1}\cdot_{\alpha}e_{2}=ae_{1},\ \
e_{2}\cdot e_{2}=e_{1}+e_{2},\ \ e_{3}\cdot e_{2}=\mu e_{3},\ \ [e_{1},
e_{2}]=a^2e_{1}.$$
These Hom-Poisson superalgbras are not Poisson superalgebra for
$a\neq  1$, or $\mu\neq 1.$

Let $b,c\in\C^*$, define a linear map $D:\A\rightarrow \A$ by
\begin{eqnarray}
D(e_1)=(a-1)be_1,\ \ D(e_2)=be_1,\ \ D(e_3)=ce_3.
\end{eqnarray}
It is easily to check that $D$ is an even derivation of $(\A,\d, [\d,\d])$ and
commutes with $\a$. By Theorem \ref {th2} in which $\lambda=0$, we obtain a super Hom-GD bialgebra  $(\mathcal {A},[\d,\d],\circ,\a)$ satisfying the following
nonzero products
$$e_{2}\c e_{2}=ab e_{1},\ \
e_{2}\c e_{3}=c\mu e_{3},\ \ [e_{1},
e_{2}]=a^2e_{1}.$$
\end{exam}

The last construction of super Hom-GD bialgebras
comes from a Hom-Lie superalgebra with a suitable linear selfmap, extending Yau' s construction of Hom-Novikov algebras in \cite{Yau}.
For a Hom-Lie superalgebra $(\LL, [\cdot, \cdot], \alpha)$, write $Z(\a(\LL))$ for the subset of $\LL$
consisting of elements $x$ such that $[x, \a(y)] = 0$ for all $y\in\LL$.

\begin{theo} \label{th5} Let $(\LL, [\cdot, \cdot], \alpha)$ be a Hom-Lie superalgebra with an even linear map $f: \LL\rightarrow \LL$
such that $f\a=\a f$. Define
\begin{eqnarray}\label{st1}
x\star y = [f(x),y],\ \ \  x\star' y = [x,f(y)], \ \ \forall \ x, y\in\LL.
\end{eqnarray}
Then
\begin{itemize}
\item[\rm(1)]  $(\LL,[\cdot, \cdot],\star,\a)$ is a super Hom-GD bialgebra if and only if the following conditions hold for all homogeneous elements $x, y, z \in \LL$:
\begin{eqnarray}
&&f([f(x),y]+[x,f(y)])-[f(x),f(y)]\in Z(\a(\LL)),\label{st2}\\
&&[f([f(x),y]), \a(z)]=(-1)^{|y||z|}[f([f(x),z]),\a (y)]. \label{st3}
\end{eqnarray}
\item[\rm(2)]  $(\LL,[\cdot, \cdot],\star',\a)$ is a super Hom-GD bialgebra if and only if the following conditions hold for homogeneous elements $x, y, z \in \LL$:
\begin{eqnarray}
[[x,f(y)]+[f(x),y],f(\a(z))]-[\a(x),f([y,f(z)])]+(-1)^{|x||y|}[\a(y),f([x,f(z)])]=0,\label{st4}
\end{eqnarray}
\begin{eqnarray}
&&[f(y),f(z)]\in Z(\a(L)),\label{st5}\\
&&[y,f(z)]+[f(y),z]- f([y,z])\in Z(\a(L)).\label{st6}
\end{eqnarray}
\end{itemize}
\end{theo}
\noindent{\it Proof.~}
(1) With a direct computation we have
\begin{eqnarray*}
\mbox{\eqref{Hom-Nov1++} holds} &\iff& [f([f(x),y]),\a(z)]-[f(\a(x)),[f(y),z]]\\
&&=(-1)^{|x||y|}\big([f([f(y),x]),\a(z)]-[f(\a(y)),[f(x),z]]\big)\\
&\iff&[f([f(x),y]+[x,f(y)])-[f(x),f(y)],\a(z)]=0\\
&\iff& \mbox{\eqref{st2} holds}.
\end{eqnarray*}
Relation \eqref{st3} follows directly by expanding \eqref{Hom-Nov2++} for the multiplication $\star$. Thus $(\LL,\star,\a)$ is a Hom-Novikov superalgebra if and only if \eqref{st2} and \eqref{st3} hold.
It remains to check the compatibility condition \eqref{HGB}. By \eqref{st3} and \eqref{super-JJ}, we have
\begin{eqnarray*}
&&[x\star y,\a(z)]-(-1)^{|y||z|}[x\star z,\a(y)]+[x,y]\star \a(z)-(-1)^{|y||z|}[x,z]\star \a(y)-\a(x)\star[y,z]\\
&=&[[f(x),y],\a(z)]-(-1)^{|y||z|}[[f(x),z],\a(y)]+[f([x,y]),\a(z)]-(-1)^{|y||z|}[f([x,z]),\a(y)]\\
&&-[f(\a(x)),[y,z]]\\
&=&[[f(x),y],\a(z)]-(-1)^{|y||z|}[[f(x),z],\a(y)]+(-1)^{|y||z|}[f([x,z]),\a(y)]\\
&&-(-1)^{|y||z|}[f([x,z]),\a(y)]-[f(\a(x)),[y,z]] \\
&=&[[f(x),y],\a(z)]+(-1)^{|x||y|}[\a(y),[f(x),z]]-[\a(f(x)),[y,z]]=0.
\end{eqnarray*}
This proves \eqref{HGB} and the first assertion.

(2) Relation \eqref{st4} is nothing but expansion of \eqref{Hom-Nov1++} for the multiplication $\star'$. By \eqref{st1} and the Hom-Jacobi identity \eqref{super-JJ},
\begin{eqnarray*}
\mbox{\eqref{Hom-Nov2++} holds}
&\iff& [[x,f(y)], f(\a(z))]=(-1)^{|y||z|}[[x,f(z)],f(\a (y))]\\&\iff& [\a(x),[f(y),f(z)]]=0\\&\iff& \mbox{\eqref{st5} holds}.
\end{eqnarray*}
Hence $(\LL,\star',\a)$ is a Hom-Novikov superalgebra if and only if \eqref{st4} and \eqref{st5} hold. Using \eqref{super-JJ} again, we have
\begin{eqnarray*}
&&[x\star' y,\a(z)]-(-1)^{|y||z|}[x\star' z,\a(y)]+[x,y]\star'\a(z)-(-1)^{|y||z|}[x,z]\star'\a(y)-\a(x)\star'[y,z]\\
&=&[[x,f(y)],\a(z)]-(-1)^{|y||z|}[[x,f(z)],\a(y)]+[[x,y],f(\a(z))]\\
&&-(-1)^{|y||z|}[[x,z],f(\a(y))]-[\a(x),f([y,z])]\\
&=&[\a(x),[y,f(z)]]+[\a(x),[f(y),z]]-[\a(x),f([y,z])]\\
&=&[\a(x),[y,f(z)]+[f(y),z]-f([y,z])].
\end{eqnarray*}
Thus \eqref{HGB} holds for the
multiplication $\star'$ if and only if
\begin{eqnarray*}
[y,f(z)]+[f(y),z]-f([y,z])\in Z(\a(L)),
\end{eqnarray*}
which is exactly \eqref{st6}. This proves the second assertion of the theorem. \QED\vskip8pt

A {\it derivation} on a Hom-Lie superalgebra is defined in the usual way. The following result follows immediately from
Theorem \ref{th5}(2).

\begin{coro}Let $(\LL, [\cdot, \cdot], \alpha)$ be a Hom-Lie superalgebra with an even derivation $d$ such that $d\a=\a d$. Define
\begin{eqnarray}\label{st7}
x\star' y = [x,d(y)], \ \forall  \ x,\ y\in\LL.
\end{eqnarray}
Then $(\LL,[\cdot, \cdot],\star',\a)$ is a super Hom-GD bialgebra if and only if the following conditions
\begin{eqnarray*}
&&[d([x,y]),d(\a(z))]-[\a(x),d([y,f(z)])]+(-1)^{|x||y|}[\a(y),d([x,d(z)])]=0,\\
&&\ \ \ \ \ \ \ \ \ \ \ \ \ \ \ \ \ \  \quad  [d(x),d(y)]\in Z(\a(L)), \ \forall \ x, y, z \in \LL.
\end{eqnarray*}
\end{coro}

In the following we assume that $\mathcal {A}=\A_0\oplus \A_1$ is a superalgebra endowed with two bilinear operations $[\d,\d]$ and
$\circ$, and an even linear selfmap $\a$. Let
\begin{eqnarray*}
\mathcal {L}(\A)=\A_0\otimes\mathbb{C}[t,t^{-1}]\oplus \A_1\otimes t^{\frac12}\mathbb{C}[t,t^{-1}]
\end{eqnarray*}
be the affinization of $\A$ equipped with a bilinear operation $[-,-]: \mathcal {L}(\A)\times\mathcal {L}(\A)\rightarrow\mathcal {L}(\A)$ defined by ($m,n\in\mathbb{Z}$ and $u,\ v\in\A$)
\begin{eqnarray}\label{aff-aff}
[u\otimes t^{m+\frac{|u|}{2}}, v\otimes t^{n+\frac{|v|}{2}}]&=&[u,v]\o t^{m+n+\frac{|u|+|v|}{2}}+(m+\frac{|u|}{2}) u\c v\o t^{m+n+\frac{|u|+|v|}{2}-1}\nonumber\\&&-(-1)^{|u||v|}(n+\frac{|v|}{2})
v\c u\o t^{m+n+\frac{|u|+|v|}{2}-1}.
\end{eqnarray}
For $a\in\A_0$, $u\in\A_1$, $m\in\mathbb{Z}$, denote
\begin{eqnarray}
a[m]=a\otimes t^{m}, \ \ \ u[m]=u\otimes t^{m+\frac12}.
\end{eqnarray}
Relation \eqref{aff-aff} is rewritten as ($a,b\in\A_0$, $u,v\in\A_1$, $m,n\in\mathbb{Z}$)
\begin{eqnarray}
[a[m],b[n]]&=&[a,b][m+n]+(m a\c b-n b\c a)[m+n-1],\label{aff1}\\
{[a[m],u[n]]}&=&[a,u][m+n]+(m a\c u-(n+\frac12) u\c a)[m+n-1],\label{aff2}\\
{[u[m],v[n]]}&=&[u,v][m+n+1]+((m+\frac12) u\c v+(n+\frac12) v\c u)[m+n].\label{aff3}
\end{eqnarray}
 Furthermore, set
\begin{eqnarray}\label{ge}
a(z)=\mbox{$\sum_{m\in\mathbb{Z}}$}a[m]z^{-m-1},\ \ u(z)=\mbox{$\sum_{n\in\mathbb{Z}}$}u[n]z^{-n-\frac32}.
\end{eqnarray}
Then \eqref{aff1}--\eqref{aff3} are respectively equivalent to
\begin{eqnarray*}\label{dis}
[a(z_1),b(z_2)]&=&[a, b](z_2)\delta(z_1-z_2)+\big(\partial_{z_2}(b\c
a)(z_2)\big)\delta(z_1- z_2)\nonumber\\&&+(a\c b+b\c
a)(z_2)\big(\partial_{z_2}\delta(z_1-z_2)\big),\\
{[a(z_1),u(z_2)]}&=&[a, u](z_2)\delta(z_1-z_2)+\big(\partial_{z_2}(u\c
a)(z_2)\big)\delta(z_1- z_2)\nonumber\\&&+(a\c u+u\c
a)(z_2)\big(\partial_{z_2}\delta(z_1-z_2)\big),\\
{[u(z_1),v(z_2)]}&=&[u, v](z_2)\delta(z_1-z_2)(\frac{z_1}{z_2})^{-\frac12}-\big(\partial_{z_2}(v\c u)(z_2)\big)\delta(z_1- z_2)(\frac{z_1}{z_2})^{-\frac12}\nonumber\\
&&+(u\c v-v\c u)(z_2)\delta(z_1-z_2)\big(\partial_{z_2}(\frac{z_1}{z_2})^{-\frac12}\big)\nonumber\\
&&+(u\c v-v\c
u)(z_2)\big(\partial_{z_2}\delta(z_1-z_2)\big)(\frac{z_1}{z_2})^{-\frac12},
\end{eqnarray*}
where $\delta(z_1-z_2)=\mbox{$\sum_{n\in\Z}$} z_1^{-n-1} z_2^{n}.$ Define an even linear map $\varphi:\mathcal {L}(\A)\rightarrow \mathcal {L}(\A)$ by
\begin{eqnarray}\label{map}
\varphi(u\otimes t^{m+\frac{|u|}{2}})=\a(u)\otimes t^{m+\frac{|u|}{2}}, \ \forall \ u\in \A,\ m\in\mathbb{Z}.
\end{eqnarray}

 The following result is a super version of \cite[Theorem 5.1]{Yuan2}.
\begin{theo}\label{h-h} $(\L,[-,-],\varphi)$
is a Hom-Lie superalgebra if and only if $(\mathcal
{A},[\d,\d],\circ,\a)$ is a  super Hom-GD bialgebra.
\end{theo}
\noindent{\it Proof.~} Assume that
$(\mathcal {A},[\d,\d],\circ,\a)$ is a  super Hom-GD bialgebra. For all $u, v, w\in\A$, $m,n,k\in\mathbb{Z}$,
\begin{eqnarray*}
[u\otimes t^{m+\frac{|u|}{2}}, v\otimes t^{n+\frac{|v|}{2}}]
&=&-(-1)^{|u||v|}[v,u]\o t^{m+n+\frac{|u|+|v|}{2}}+(m+\frac{|u|}{2}) u\c v\o t^{m+n+\frac{|u|+|v|}{2}-1}\nonumber\\&&-(-1)^{|u||v|}(n+\frac{|v|}{2})
v\c u\o t^{m+n+\frac{|u|+|v|}{2}-1}\\
&=&-(-1)^{|u||v|}[v\otimes t^{n+\frac{|v|}{2}}, u\otimes t^{m+\frac{|u|}{2}}]
\end{eqnarray*}
by \eqref{aff-aff}. Thus \eqref{super-ss} holds. To check \eqref{super-JJ}, we compute
\begin{eqnarray}\label{jj}
&&(-1)^{|u||w|}[[u\otimes t^{m+\frac{|u|}{2}}, v\otimes t^{n+\frac{|v|}{2}}],\varphi(w\o t^{k+\frac{|w|}{2}})]\nonumber\\
&& +(-1)^{|v||u|}[[v\o t^{n+\frac{|v|}{2}},w\o
t^{k+\frac{|w|}{2}}],\varphi(u\o t^{m+\frac{|u|}{2}})]\nonumber\\
&&+(-1)^{|v||w|}[[w\o t^{k+\frac{|w|}{2}},u\o t^{m+\frac{|u|}{2}}],\varphi(v\o
t^{n+\frac{|v|}{2}})]\nonumber\\
&=&\D_1\o t^{m+n+k+\frac{|u|+|v|+|w|}{2}}+\D_2\o t^{m+n+k+\frac{|u|+|v|+|w|}{2}-1}
+\D_3\o
t^{m+n+k+\frac{|u|+|v|+|w|}{2}-2},
\end{eqnarray}
where
\begin{eqnarray}
\D_1&=&(-1)^{|u||w|}[[u,v],\a(w)]+(-1)^{|u||v|}[[v,w],\a(u)]+(-1)^{|v||w|}[[w,u],\a(v)],\label{de1}\\
\D_2&=&(-1)^{|u||w|}(m+\frac{|u|}{2})\big([u,v]\c\a(w)+[u\c
v,\a(w)]-\a(u)\c[v,w]\nonumber\\&&\ \ \ \ \ \ \ \ \ \ \ \ \ -(-1)^{|v||w|}[u,w]\c\a(v)-(-1)^{|v||w|}[u\c w,\a(v)]\big)\nonumber\\
&&-(-1)^{|u|(|w|+|v|)}(n+\frac{|v|}{2})\big([v,u]\c\a(w)+[v\c u,\a(w)]-\a(v)\c[u,w]\nonumber\\&&\ \ \ \ \ \ \ \ \ \ \ \ \ -(-1)^{|u||w|}[v,w]\c\a(u)-(-1)^{|u||w|}[v\c
w,\a(u)]\big)\nonumber\\
&&-(-1)^{|v||w|}(k+\frac{|w|}{2})\big(\a(w)\c[u,v]+(-1)^{|u||v|}[w,v]\c\a(u)+(-1)^{|u||v|}[w\c
v,\a(u)]\nonumber\\&&\ \ \ \ \ \ \ \ \ \ \ \ \ -[w,u]\c\a(v)-[w\c u,\a(v)]\big),\label{de2}\\
\D_3&=&(-1)^{|u||w|}[(m+\frac{|u|}{2})^2-(m+\frac{|u|}{2})]\big((u\c v)\c\a(w)-(-1)^{|v||w|}(u\c
w)\c\a(v)\big)\nonumber\\&&+(-1)^{|u||v|}[(n+\frac{|v|}{2})^2-(n+\frac{|v|}{2})]\big((v\c w)\c\a(u)-(-1)^{|u||w|}(v\c
u)\c\a(w)\big)\nonumber\\
&&+(-1)^{|v||w|}[(k+\frac{|w|}{2})^2-(k+\frac{|w|}{2})]\big((w\c u)\c\a(v)-(-1)^{|u||v|}(w\c
v)\c\a(u)\big)\nonumber\\
&&+(m+\frac{|u|}{2})(n+\frac{|v|}{2})\big((-1)^{|u||w|}(u\c v)\c\a(w)-(-1)^{|u|(|v|+|w|)}(v\c
u)\c\a(w)\nonumber\\&&\ \ \ \ \ \ \ \ \ \ \ \ \ \ \ \ \ \ \ \ \ \ \ -(-1)^{|u||w|}\a(u)\c(v\c w)+(-1)^{|u|(|v|+|w|)}\a(v)\c(u\c w)\big)\nonumber\\
&&+(m+\frac{|u|}{2})(k+\frac{|w|}{2})\big((-1)^{|v||w|}(w\c
u)\c\a(v)-(-1)^{|w|(|u|+|v|)}(u\c w)\c\a(v)\nonumber\\&&\ \ \ \ \ \ \ \ \ \ \ \ \ \ \ \  \ \ \ \ \ \ \ -(-1)^{|v||w|}\a(w)\c(u\c v)+(-1)^{|w|(|u|+|v|)}\a(u)\c(w\c
v)\big)\nonumber\\&&+(n+\frac{|v|}{2})(k+\frac{|w|}{2})\big((-1)^{|u||v|}(v\c w)\c\a(u)-(-1)^{|v|(|u|+|w|)}(w\c
v)\c\a(u)\nonumber\\&&\ \ \ \ \ \ \ \ \ \ \ \ \ \ \ \ \ \ \ \ \ \  -(-1)^{|u||v|}\a(v)\c(w\c u)+(-1)^{|v|(|u|+|w|)}\a(w)\c(v\c u)\big).\label{de3}
\end{eqnarray}
By the hypothesis that $(\mathcal {A},[\d,\d],\circ,\a)$ is a  super Hom-GD bialgebra, $\D_1=\D_2=\D_3=0$ and so $\eqref{jj}$ equals to zero.
Thus $(\L,[-,-],\varphi)$ is a Hom-Lie superalgebra.

If $(\L,[-,-],\varphi)$
is a Hom-Lie superalgebra, then $\eqref{jj}$ equals to zero. It follows that
 $\D_1=0$ (so $(\mathcal
{A}, [\d,\d],\a)$ is a Hom-Lie superalgebra), and $\D_2=\D_3=0$ for all $m,n,k\in\mathbb{Z}$.
Then \eqref{de3} implies \eqref{Hom-Nov1++} and \eqref{Hom-Nov2++}, whereas  \eqref{de2} gives \eqref{HGB}.
Thus $(\mathcal
{A},[\d,\d],\circ,\a)$ is a super Hom-GD bialgebra.  \QED\vskip8pt

\begin{exam}\label{exam4-2}\rm By Example \ref{exam2}, there exsits a super Hom-GD bialgebra $(\A,[\cdot,\cdot],\circ,\alpha)$ with $\A_0=\C x_1\oplus \C x_2$, $\A_1=\C y$, satisfying the following nonzero products
\begin{eqnarray*}
x_1\circ x_1=x_2, \ \ x_2\circ x_2=x_1, \ \ x_1\circ y=y,\ \
 x_2\circ y=y,\ \ [x_1,y]=y, \ [x_2,y]=y,
\end{eqnarray*}
and $\alpha(x_1)=x_2,\  \alpha(x_2)=x_1,\ \alpha(y)=0$.
Let $\mathcal {L}(\A)$ be the affinization of $\A$. By \eqref{aff1}--\eqref{aff3}, the following nontrivial relations hold for $m,n,k\in\Z$:
\begin{eqnarray*}
&&[x_1[m],x_1[n]]=(m-n)x_2[m+n-1],\ \ [x_1[m],y[k]]=y[m+k]+my[m+k-1],\\
&&{[x_2[m],x_2[n]]}=(m-n)x_1[m+n-1],\ \ {[x_2[n],y[k]]}=y[n+k]+ny[n+k-1].
\end{eqnarray*}
By Theorem \ref{h-h}, we obtain a Hom-Lie superalgebra $(\mathcal {L}(\A), [\cdot,\cdot],\varphi)$, where
$\varphi$ is defined by
\begin{eqnarray*}
\varphi(x_1[m])=x_2[m], \ \varphi(x_2[n])=x_1[n],\ \varphi(y[k])=0.
\end{eqnarray*}
\end{exam}
\begin{exam}\rm Consider the super Hom-GD bialgebra $(\A,[\cdot,\cdot],\circ,\alpha)$ from Example \ref{exam3}.
According to Theorem \ref{h-h}, we obtain a Hom-Lie superalgebra $(\mathcal {L}(\A), [\cdot,\cdot],\varphi)$ satisfying the following nonzero products
\begin{eqnarray*}
[x_1[m],x_2[n]]&=&\lambda^2x_1[m+n]+\frac{m+n}{2}\lambda^2x_1[m+n-1],\\
{[x_2[m],x_2[n]]}&=&\frac{m-n}{2}x_2[m+n-1],\\
{[x_2[n],y[k]]}&=&-\frac\lambda 2 y[n+k]-(k+\frac{1}{2})\frac\lambda 2 y[n+k-1],\\
{[y[m],y[k]]}&=&(m+k+1)\frac{\lambda^2}{2} x_1[m+k],
\end{eqnarray*}
and $\varphi(x_1[m])=\lambda^2 x_1[m], \ \varphi(x_2[n])=x_2[n], \ \varphi(y[k])=\l y[k],$ for $m,n,k\in\Z.$
\end{exam}

\vs{8pt}\

\cl{\bf\S4. \ Hom-Lie conformal superalgebras
}\setcounter{section}{4}\setcounter{equation}{0}\setcounter{theo}{0}
\vskip8pt

In this section, we introduce the notion of Hom-Lie conformal superalgebra, which is a Hom-generalization of Lie conformal superalgebras and also a superanalogue of Hom-Lie conformal algebras. In particular, we introduce quadratic Hom-Lie conformal superalgebras and establish equivalence of quadratic Hom-Lie conformal superalgebras and super Hom-GD bialgebras.

\begin{defi}\label{HLCSA}{\rm
A Hom-Lie conformal superalgebra $\mathcal
{R}=\mathcal {R}_0\oplus \mathcal {R}_1$ is a $\C[\partial ]$-module equipped with an even linear selfmap $\a$ satisfying $\a\p=\p\a$, and a $\C$-bilinear map
 \begin{eqnarray*}
 \mathcal {R}\otimes \mathcal {R}\rightarrow \C[\lambda]\otimes \mathcal {R}, \ \ a\otimes b\mapsto[a_\l b],
 \end{eqnarray*}
such that the following axioms hold for all homogenous elements $a, b, c\in \mathcal {R}$:
\begin{eqnarray}
[\partial a_\lambda b]&=&-\lambda[a_\lambda b],\ [a_{\l}\p b]=(\p+\l)[a_{\l} b],\ \ \mbox{(conformal\  sesquilinearity)}\label{HL1}\\
{[a_\lambda b]} &=& -(-1)^{|a||b|}[b_{-\lambda-\partial}a],\ \ \ \ \ \ \mbox{(skew-symmetry)}\label{HL2}\\
{[\a(a)_\lambda[b_\mu c]]}&=&[[a_\lambda b]_{\lambda+\mu
}\a(c)]+(-1)^{|a||b|}[\a(b)_\mu[a_\lambda c]].\ \ \mbox{(Hom-Jacobi \
identity)}\label{HL3}
\end{eqnarray}
}\end{defi}

  A Hom-Lie
conformal superalgebra $\mathcal {R}$ is called {\it finite} if
$\mathcal {R}$ is a finitely generated $\C[\partial]$-module. The
{\it rank} of $\mathcal {R}$ is its rank as a $\C[\partial]$-module.
\begin{rema}{\rm We recover Lie conformal superalgebras when $\a={\rm id}.$ The
Hom-Lie conformal algebras are obtained when the odd part is trivial. In addition, we can associate to any Lie conformal superalgebra a Hom-Lie conformal superalgebra structure by
taking $\a=0$.
}\end{rema}

The following notion is a Hom-analogue of quadratic conformal superalgebra studied in \cite{X1}.
\begin{defi}{\rm Let
($\mathcal {R},[\cdot_\lambda\cdot],\a$) be a Hom-Lie conformal superalgebra. If $\mathcal {R}=\C[\p]V$ is a free
 $\C[\partial ]$-module over a superspace $V$ and the $\l$-bracket is of the following form
\begin{eqnarray}
[a_\l b]=\p u+\l v+w, \ \ \mbox{for} \ u,v,w\in V,
\end{eqnarray} then ($\mathcal {R},[\cdot_\lambda\cdot],\a$) is called a quadratic Hom-Lie conformal superalgebra.
}\end{defi}

 We in \cite{Yuan2} gave a construction of free Hom-Lie conformal algebras from formal distribution Hom-Lie algebras, generalizing the classical construction of free Lie conformal algebras from formal distribution Lie algebras due to Kac \cite{ka}. In the following we extend this construction to Hom-Lie conformal superalgebras.

 Let $(\LL,[\d,\d],\a)$ be a Hom-Lie superalgebra. Consider the following $\LL$-valued formal distributions
 \begin{eqnarray*}
 a(z)=\mbox{$\sum\limits_{m\in\Z}$}a_mz^{-m-1}\in
\LL[[z,z^{-1}]],\ \
 b(w)=\mbox{$\sum\limits_{n\in\Z}$}b_nw^{-n-1}\in
\LL[[w,w^{-1}]].
\end{eqnarray*}
Define a formal distribution in two variables by
\begin{eqnarray}\label{LB2}
[a(z),b(w)]=
\mbox{$\sum\limits_{m,n\in\Z}$}[a_{m},b_{n}]z^{-m-1}w^{-n-1}\in
\LL[[z,z^{-1},w,w^{-1}]].
\end{eqnarray}
If there exists a positive integer $N$ such that $(z-w)^N [a(z),b(w)] = 0,$  then $a(z)$ and $b(w)$ are said to be {\it mutually local}. Taking the formal Fourier transform (cf. \cite{ka})
\begin{eqnarray*}
F^\lambda_{z,w}={\rm Res}_ze^{\lambda(z-w)}={\rm Res}_z\mbox{$\sum_{j\in \Z^+
}$}\frac{\lambda^j}{j!}(z-w)^j
\end{eqnarray*}
on both parts of \eqref{LB2}, we obtain a $\mathbb{C}$-bilinear map
\begin{eqnarray*}
[\d\,_{\lambda}\d]:\LL[[w,w^{-1}]]\otimes
\LL[[w,w^{-1}]]\rightarrow\LL[[w,w^{-1}]][[\lambda]]\end{eqnarray*} such that
\begin{eqnarray}\label{lamda}
[a(w)_\lambda b(w)]= F^\lambda_{z,w}[a(z),b(w)],
\end{eqnarray}
which is called a {\it $\lambda$-bracket} between $a(w)$ and $b(w)$, and simply denoted by
$[a_\lambda b]$. If $a(z)$ and $b(w)$ are mutually local, then Kac's decomposition theorem (cf. \cite{ka}) gives
\begin{eqnarray}\label{sum+}
[a(z),b(w)]=\mbox{$\sum\limits_{j\in\Z^+}$}a(w)_{(j)}b(w)\frac{\p_w^{j}\delta(z,w)}{j!},
\end{eqnarray}
where
\begin{eqnarray}\label{jsum+}
a(w)_{(j)}b(w)={\rm Res}_z(z-w)^j[a(z),b(w)]
\end{eqnarray}
is called a {\it $j$-product} of $a(w)$ and $b(w)$, and simply denoted by $a_{(j)}b$. By \eqref{lamda} and \eqref{jsum+}, the $\lambda$-bracket
is related to the $j$-products as follows:
\begin{eqnarray}\label{lj}
[a_\lambda
b]=\mbox{$\sum\limits_{j\in\Z^+}$}\frac{\lambda^j}{j!}(a_{(j)}b).
\end{eqnarray}
Furthermore, define actions of $\p$ and $\a$ on
$a(z)=\mbox{$\sum_{m\in\Z}$}a_{m}z^{-m-1}\in\LL[[z,z^{-1}]]$ by
\begin{eqnarray}
(\p a)(z)=\p_z(a(z)), \ \ \a(a)(z)=\mbox{$\sum_{m\in\Z}$}\a(a_{m})z^{-m-1}.
\end{eqnarray}

\begin{prop} \label{p*} The $\l$-bracket satisfies the following properties:
\begin{itemize}
\item[\rm (1)]$[\partial a_\lambda
b]=-\lambda[a_\lambda b],\ \ {[a_\lambda
\partial b]}=(\p+\lambda)[a_\lambda b].$
\item[\rm(2)] if $a(z)$ and $b(w)$ are mutually local, then ${[a_\lambda b]}=-(-1)^{|a||b|}[b_{-\lambda-\partial}a].$
\item[\rm(3)] ${[\a(a)_\lambda[b_\mu c]]}=[[a_\lambda b]_{\lambda+\mu
}\a(c)]+(-1)^{|a||b|}[\a(b)_\mu[a_\lambda c]].$
\end{itemize}
\end{prop}
\noindent{\it Proof.}\ \ It can be similarly proved to \cite[Proposition
4.6]{Yuan2}.\QED\vskip4pt\par

A subset $F\subset \LL[[z,z^{-1}]]$ is called a {\it local family}
of $\LL$-valued formal distributions if all pairs of its
constituents are mutually local.
\begin{defi}{\rm Let $(\LL,[\d,\d],\a)$ be a Hom-Lie superalgebra. If there exists a local family $F$ of $\LL$-valued formal
distributions with their Fourier coefficients generating the whole
$\LL$, then $F$ is said to endow $\LL$ with a structure of
formal distribution Hom-Lie superalgebra. In this case, we denote
by $(\LL,F)$ to emphasize the role of $F$.}
\end{defi}

\begin{defi} {\rm Let $(\LL, [\d,\d],\a)$ be a Hom-Lie superalgebra. A local family
$F\subset \LL[[z,z^{-1}]]$ is called a
conformal family if it is closed under their $j$-products and invariant under the actions of $\p$ and $\a$.}
\end{defi}

It was pointed out in \cite{Yuan2} that if $(\LL, F)$ is a formal
distribution Hom-Lie algebra, one can always include $F$ in the minimal conformal family $\bar {F}$, such that
$(\LL, F)$ can be endowed with a free Lie conformal algebra structure $\mathcal{R}=\C[\p]\bar F,$ such that
$[a_\l b]=F^{\l}_{z,w}[a(z),b(w)]$, $\p=\p_z$ and $\a(f(\p)u)=f(\p)\a(u)$, for $f(\p)\in\C[\p],$ $u\in \bar F$.
Such a construction can be naturally extended to formal
distribution Hom-Lie superalgebras.

\begin{exam}{\rm  Let $(\LL, [\d,\d],\a)$ be a Hom-Lie superalgebra. The loop algebra
associated to $\LL$ is
$$\hat \LL=\LL\otimes \C[t,t^{-1}],\ \ \ |u\otimes t^m|=|u|, \ \ \mbox{for}\ u\in\LL,\ m\in\Z$$
with
\begin{eqnarray*}
[u\otimes t^m, v\otimes t^n]=[u,v]\otimes t^{m+n},\ \ \mbox{for} \
u, v\in \LL,\ m,n\in\Z.
\end{eqnarray*}
Extend $\a$ to $\hat \LL$ by $\a(u\otimes t^m)=a(u)\otimes t^m$.
Then $(\hat \LL, [\d,\d],\a)$ is a Hom-Lie superalgebra.
We introduce the family $F$ of $\hat \LL$-valued formal distributions
$$u(z)=\mbox{$\sum\limits_{n\in\Z}$}(u\otimes t^n)z^{-n-1}, \ \mbox{for}\ u\in\LL.$$
It is easily verified that
\begin{eqnarray*}
[u(z),v(w)]=[u,v](w)\delta(z,w),
\end{eqnarray*}
which is equivalent to $u_{(0)}v=[u,v]\in \LL$, and $u_{(n)}v=0$ for $n>0$. Hence $(\hat \LL, F)$ is a formal distribution Hom-Lie superalgebra. The associated Hom-Lie conformal superalgebra is $\mathcal{R}=\C[\p]\LL$, with $\lambda$-bracket of the following form 
\begin{eqnarray}\label{1111}
[u_{\l}v]=[u,v], \ \mbox{for} \ u,v\in \LL.
\end{eqnarray}
Indeed, we can extend the $\lambda$-bracket and $\a$ to the whole $\mathcal{R}$ by
\begin{eqnarray}
[f(\partial) u_\lambda h(\partial)
v]&=& f(-\lambda)h(\partial+\lambda)[u_\lambda v],\label{lamda22+2}\\
\a(f(\p)u)&=&f(\p)\a(u),\label{lamda22+2a}
\end{eqnarray}
for $f(\partial),h(\partial)\in \mathbb{C}[\p]$,
$u,v\in\mathcal{R}$.
Then \eqref{HL1} follows directly from (\ref{lamda22+2}), whereas (\ref{lamda22+2a}) gives $\a\p=\p\a$. It suffices to check \eqref{HL2} and \eqref{HL3} on the generators. By \eqref{1111}, axioms \eqref{HL2} and \eqref{HL3} naturally hold. Thus $\mathcal{R}=\C[\p]\LL$ is a Hom-Lie conformal superalgebra, which is viewed as a {\it current-like  Hom-Lie conformal superalgebra}. }\end{exam}

\begin{rema}\label{remark}{\rm If $\mathcal{R}$ is a free $\C[\p]$-module with $\l$-bracket
 $[\d_{\l}\d]$ defined on a $\C[\p]$-basis of $\mathcal{R}$ such that \eqref{HL2}
and \eqref{HL3} hold, there is a unique extension of this
$\l$-bracket via \eqref{HL1} to the whole
$\mathcal{R}$ (as shown in \eqref{lamda22+2}). It is easy to see that
\eqref{HL2} and \eqref{HL3} also hold for this extension. Thus,
$\mathcal{R}$ is equipped with a Hom-Lie conformal superalgebra
structure. In the sequel, we shall often describe Hom-Lie conformal superalgebras on free $\C[\p]$-modules by defining
$\l$-bracket on a fixed $\C[\p]$-basis. }\end{rema}

\begin{exam}\label{exm4-1}{\rm
Recall that the classical super Virasoro algebra $sVir$ is
an infinite-demensional superalgebra generated by $\{L_n, F_n,
G_n|n\in Z\}$ (see \cite{MM}), such that
\begin{eqnarray}
&&[L_m,L_n]=(m-n)L_{m+n},\ \ \
[L_m,F_n]=-nF_{n+m},\ \ \
[F_m,F_n]=0,\label{sVir1}\\
&&[L_m,G_n]=(m-n)G_{n+m},\ \ \
[F_m,G_n]=G_{n+m},\ \ \ \ \ \,
[G_m,G_n]=0.\label{sVir2}
\end{eqnarray}
The $\Z_2$-grading is
defined by requiring that ${\rm deg}(L_i)={\rm deg}(F_i)=0$ and ${\rm deg}(G_i)=1$ for $i\in\Z$.
Set
\begin{eqnarray*}
L(z)=\mbox{$\sum\limits_{n\in\Z}$}L_nz^{-n-2},\ \ \
F(z)=\mbox{$\sum\limits_{n\in\Z}$}F_nz^{-n-1},\ \ \
G(z)=\mbox{$\sum\limits_{n\in\Z}$}G_nz^{-n-2}.
\end{eqnarray*}
By a straightforward comuputation, relations \eqref{sVir1} and \eqref{sVir2} amount to the following
$\l$-brackets
\begin{eqnarray}
&&[L_\lambda L]=(\p+2\lambda) L,\ \ [L_\lambda G]=(\p+2\lambda)G,\ \ [G_\lambda L]=(\p+2\lambda)G,\ \ [F_\lambda G]=G,\label{a21}\\
&&{[L_\lambda F]}=(\p+\lambda)F,\ \ \ [F_\lambda L]=\lambda F,\ \ \ [F_\lambda F]=0,\ \ \ [G_\lambda G]=0,
 \ \ \ [G_\lambda F]=-G.\label{a22}
\end{eqnarray}
One can check that \eqref{a21} and  \eqref{a22} define a free Lie conformal superalgebra
$\mathfrak{C}sVir=\C[\p]L\oplus\C[\p]F\oplus\C[\p]G$.
Define an even linear map $\a:\mathfrak{C}sVir\rightarrow \mathfrak{C}sVir $ by
\begin{eqnarray}\label{alpha}
\a(L)=f_1(\p)L+f_2(\p)F, \ \
{\a(F)}=g_1(\p)L+g_2(\p)F,\ \
{\a(G)}=h(\p)G,
\end{eqnarray}
where $f_1(\p),f_2(\p),g_1(\p),g_2(\p),h(\p)\in\C[\partial]$.
In the following we aim to find restrictions on $\a$ such that $(\mathfrak{C}sVir,[\cdot_\l\cdot],\a)$ forms a Hom-Lie conformal superalgebra.

Applying the Hom-Jacobi identity to the triple $(L,L,L)$ gives
\begin{eqnarray*}
[\a(L)_\lambda[L_\mu L]]=[[L_\lambda L]_{\lambda+\mu}\a(L)]+[\a(L)_\mu[L_\lambda L]],\\
\end{eqnarray*}
where the left-hand side reads
\begin{eqnarray}\label{a1}
[\a(L)_\lambda[L_\mu L]]=(\p+\lambda+2\mu)\big(f_1(-\lambda)(\p+2\lambda)L+f_2(-\lambda)\lambda F\big),
\end{eqnarray}
and the right-hand side reads
\begin{eqnarray}\label{a2}
&&{[[L_\lambda L]_{\lambda+\mu}\a(L)]}+{[\a(L)_\mu[L_\lambda L]]}\nonumber\\
&&\ \ \ \ \ =(\lambda-\mu)\big(f_1(\p+\lambda+\mu)(\p+2\lambda+2\mu)L+
f_2(\p+\lambda+\mu)(\p+\lambda+\mu)F\big)\nonumber\\&&
\ \ \ \ \ \ \ \ +(\p+\mu+2\lambda)\big(f_1(-\mu)(\p+2\mu)L+f_2(-\mu)\mu F\big).
\end{eqnarray}
Combining \eqref{a1} with \eqref{a2} yields
\begin{eqnarray}
f_1(-\lambda)(\p+\lambda+2\mu)(\p+2\lambda)&=&f_1(\p+\lambda+\mu)(\lambda-\mu)(\p+2\lambda+2\mu)\nonumber\\
&&+f_1(-\mu)(\p+\mu+2\lambda)(\p+2\mu),\label{a3}\\
{f_2(-\lambda)(\p+\lambda+2\mu)\lambda}&=&f_2(\p+\lambda+\mu)(\lambda-\mu)(\p+\lambda+\mu)\nonumber\\
&&+f_2(-\mu)(\p+\mu+2\lambda)\mu.\label{a4}
\end{eqnarray}
Setting $\p=-2\mu$ in \eqref{a3}, we have
\begin{eqnarray*}
f_1(-\lambda)=f_1(\lambda-\mu),
\end{eqnarray*}
which implies $f_1(\lambda)=a$, for some $a\in\C$. Taking $\mu=0$ in \eqref{a4} gives
\begin{eqnarray*}
f_2(-\lambda)=f_2(\p+\lambda),
\end{eqnarray*}
and thus $f_2(\p)=b$ for some $b\in\C$. By \eqref{alpha}, we get
\begin{eqnarray}\label{L}
\a(L)=aL+bF, \ \mbox{for\ some} \ a, b\in\C.
\end{eqnarray}
Applying the Hom-Jacobi identity to $(L,L,F)$, together with \eqref{L}, gives
\begin{eqnarray}
[\a(L)_\lambda[L_\mu F]]&=&a(\p+\mu+\lambda)(\p+\lambda)F\nonumber\\
&=&[[L_\lambda L]_{\lambda+\mu}\a(F)]+[\a(L)_\mu[L_\lambda F]]\nonumber\\
&=&(\lambda-\mu)(g_1(\p+\mu+\lambda)(\p+2\lambda+2\mu)L
+g_2(\p+\mu+\lambda)(\p+\mu+\lambda)F)\nonumber\\
&&+a(\p+\mu+\lambda)(\p+\mu)F.\label{a5}
\end{eqnarray}
Comparing the coefficients of $L$ in \eqref{a5} gives
\begin{eqnarray*}
(\lambda-\mu)g_1(\p+\mu+\lambda)(\p+2\lambda+2\mu)=0,
\end{eqnarray*}
from which we deduce that $g_1(\p)=0$. Equating the coefficients of $F$ in \eqref{a5} gives
\begin{eqnarray}\label{a6}
a(\p+\mu+\lambda)(\p+\lambda)=(\lambda-\mu)g_2(\p+\mu+\lambda)(\p+\mu+\lambda)+a(\p+\mu+\lambda)(\p+\mu),
\end{eqnarray}
from which we obtain $g_2(\lambda)=a$ by setting $\p=-\mu$. Finally, applying the Hom-Jacobi identity to the triple $(L,F,G)$ gives
\begin{eqnarray*}
[\a(L)_\lambda[F_\mu G]]&=&a(\p+2\lambda)G+bG\\
&=&{[[L_\lambda F]_{\lambda+\mu}\a(G)]}+{[\a(F)_\mu[L_\lambda G]]}\\
&=&-\mu h(\p+\mu+\lambda)G+a(\p+2\lambda+\mu)G.
\end{eqnarray*}
It follows that $a(\p+2\lambda)+b=-\mu h(\p+\mu+\lambda)+a(2\lambda+\p+\mu)$, which gives
$h(\p)=a$, and $b=0$.

To sum up, $(\mathfrak{C}sVir,[\cdot_\l\cdot],\a)$ is a Hom-Lie conformal superalgebra if and only if
the even linear map $\a:\mathfrak{C}sVir\rightarrow \mathfrak{C}sVir $ satisfies
\begin{eqnarray*}\label{alpha111}
\a(L)=aL, \ \ {\a(F)}=aF,\ \ {\a(G)}=aG, \ \mbox{for\ some}\ a\in\C.
\end{eqnarray*}
}\end{exam}

The following result is a superanalogue of \cite[Theorem 5.2]{Yuan2}, generalizing Gel'fand and Dorfman's statement on equivalence of Gel'fand-Dorfman bialgebras and Lie conformal algebras.
\begin{theo}\label{main} Let $\mathcal {A}$ be a superspace equipped with two operations $\circ$ and $[\d,\d]$, and an even linear selfmap $\a$. Let $\mathcal{R}=\mathbb{C}[\p]\A$ be the free $\C[\p]$-module over $\A$. Extend $\a$ to $\mathcal{R}$ by
\begin{eqnarray}
\label{*} \a (f(\p) u)=f(\p)\a(u), \ \forall\ f(\p)\in\C[\p], \ u\in\A,
\end{eqnarray}
and define $\l$-bracket $[\d_\l\d]:\A\otimes\A\rightarrow\A[\l]$ by
\begin{eqnarray}\label{**}[u_{\lambda}v]=[v,u]+(\partial+\lambda)(v\c
u)+(-1)^{|u||v|}\lambda(u\c v), \ \forall \ u, v\in\A.
\end{eqnarray}
Then $(\mathcal{R},[\d_\l\d],\a)$ is a quadratic
Hom-Lie conformal superalgebra if and only if $(\mathcal{A},[\d,\d],\circ,\a)$ is a super Hom-GD bialgebra.
\end{theo}
\noindent{\it Proof.~}Suppose that $(\mathcal {A},[\d,\d],\circ,\a)$
is a super Hom-GD bialgebra. It follows from \eqref{*} that $\a\p=\p\a$. Extend the $\lambda$-bracket defined by \eqref{**} to $\mathcal{R}$
by
\begin{eqnarray*}
[f(\partial) u_\lambda h(\partial)
v]= f(-\lambda)h(\partial+\lambda)[u_\lambda v],  \ \forall \ f(\partial),h(\partial)\in \mathbb{C}[\p], \
u,v\in\mathcal{A}.
\end{eqnarray*}
Thus \eqref{HL1} is naturally satisfied. Then it
suffices to verify \eqref{HL2} and \eqref{HL3} on the generators. With a direct computation we have
\begin{eqnarray}\label{skew-s}
[v_{-\lambda-\p}u]&=&[u,v]+(-\lambda-\p+\p)(u\c v)+(-\lambda-\p)(-1)^{|u||v|}(v\c u)\nonumber\\
&=&-(-1)^{|u||v|}[v,u]-\lambda(u\c v)-(\lambda+\p)(-1)^{|u||v|}(v\c u)\nonumber\\
&=&-(-1)^{|u||v|}[u_\l v],
\end{eqnarray}
which proves \eqref{HL2}. To show \eqref{HL3}, we compute separately
\begin{eqnarray}
[\a(u)_\l [v_{\mu}
w]]
&=&[[w,v],\a(u)]+\p\big([w,v]\c\a(u)+[w\c v,\a(u)]\big)+\p^2(w\c v)\c\a(u)\nonumber\\
&&+\l\big([w,v]\c\a(u)+(-1)^{|u|(|v|+|w|)}\a(u)\c[w,v]+[w\c v,\a(u)]\big)\nonumber\\&&+\mu\big([w\c v,\a(u)]+(-1)^{|v||w|}[v\c w,\a(u)]\big)\nonumber\\
&&+\p\l\big(2(w\c v)\c\a(u)+(-1)^{|u|(|v|+|w|)}\a(u)\c(w\c v)\big)\nonumber\\&&+\l^2((w\c
v)\c\a(u)+(-1)^{|u|(|v|+|w|)}\a(u)\c(w\c v)\big)\nonumber\\
&&+\p\mu\big((w\c v)\c a(u)+(-1)^{|v||w|}(v\c w)\c a(u)\big)\nonumber\\&&+\mu\l\big((w\c
v)\c\a(u)+(-1)^{|u|(|v|+|w|)}\a(u)\c(w\c v)\nonumber\\
&&+(-1)^{|v||w|}(v\c w)\c a(u)+(-1)^{|u||v|+|u||w|+|v||w|}\a(u)\c(v\c w)\big),\label{jj1}\\
{[\a(v)_\mu [u_{\l} w]]}
&=&[[w,u],\a(v)]+\p\big([w,u]\c\a(v)+[w\c u,\a(v)]\big)+\p^2(w\c u)\c\a(v)\nonumber\\&&
+\mu\big([w,u]\c\a(v)+(-1)^{|v|(|u|+|w|)}\a(v)\c[w,u]+[w\c
u,\a(v)]\big)
\nonumber\\&&
+\l\big([w\c u,\a(v)]+(-1)^{|u||w|}[u\c w,\a(v)]\big)\nonumber\\&&+\p\mu\big(2(w\c
u)\c\a(v)+(-1)^{|v|(|u|+|w|)}\a(v)\c(w\c u)\big)
\nonumber\\&&
+\mu^2((w\c u)\c\a(v)+(-1)^{|v|(|u|+|w|)}\a(v)\c(w\c
u))\nonumber\\&&+\p\l\big((w\c u)\c\a(v)+(-1)^{|u||w|}(u\c w)\c\a(v)\big)
\nonumber\\&&
+\mu\l\big((w\c
u)\c\a(v)+(-1)^{|v|(|u|+|w|)}\a(v)\c(w\c
u)\nonumber\\&&+(-1)^{|u||w|}(u\c w)\c\a(v)+(-1)^{|u||v|+|u||w|+|v||w|}\a(v)\c(u\c
w)\big),\label{jj2}\\
{[[u_{\l} v]_{\l+\mu}\a(w)]}
&=&[\a(w),[v,u]]+\p\a(w)\c[v,u]+(-1)^{|u||v|}\p\l\a(w)\c(u\c v)\nonumber\\
&&+\mu\big(\a(w)\c[v,u]+(-1)^{|w|(|u|+|v|)}[v,u]\c\a(w)-[\a(w),v\c u]\big)\nonumber\\
&&+\l\big(\a(w)\c[v,u]+(-1)^{|w|(|u|+|v|)}[v,u]\c\a(w)+(-1)^{|u||v|}[\a(w),u\c v]\big)\nonumber\\
&&+(-1)^{|u||v|}\l^2\big(\a(w)\c(u\c v)+(-1)^{|w|(|u|+|v|)}(u\c v)\c\a(w)\big)\nonumber\\
&&-\p\mu\a(w)\c(v\c u)-\mu^2\big(\a(w)\c(v\c u)+(-1)^{|w|(|u|+|v|)}(v\c u)\c\a(w)\big)\nonumber\\
&&+\mu\l\big((-1)^{|u||v|}\a(w)\c (u\c v)+(-1)^{|u||v|+|u||w|+|v||w|}(u\c v)\c\a(w)\nonumber\\
&&-\a(w)\c(v\c u)-(-1)^{|w|(|u|+|v|)}(v\c
u)\c\a(w)\big).\label{jj3}
\end{eqnarray}
By \eqref{jj1}--\eqref{jj3}, the Hom-Jacobi identity \eqref{HL3} holds if and only if the following
relations hold
\begin{eqnarray}
0&=&[[w,v],\a(u)]-(-1)^{|u||v|}[[w,u],\a(v)]-[\a(w),[v,u]],\label{ee1}\\
0&=&{[w,v]\c\a(u)+[w\c v,\a(u)]}\nonumber\\&&-(-1)^{|u||v|}\big([w,u]\c\a(v)+[w\c
u,\a(v)]\big)-\a(w)\c[v,u],
\label{ee2}\\
0&=&{(w\c v)\c\a(u)}-(-1)^{|u||v|}(w\c u)\c\a(v),\label{ee3}
\end{eqnarray}
\begin{eqnarray}
0&=&{[w,v]\c\a(u)+(-1)^{|u|(|v|+|w|)}\big(\a(u)\c[w,v]-[u\c w,\a(v)]\big)+[w\c v,\a(u)]}-\a(w)\c[v,u]\nonumber\\&&-(-1)^{|u||v|}\big([w\c u,\a(v)]+[\a(w),u\c
v]\big)-(-1)^{|w|(|u|+|v|)}[v,u]\c\a(w), \label{ee4}\\
0&=&[w\c v,\a(u)]+(-1)^{|v||w|}[v\c w,\a(u)]-\a(w)\c[v,u]-(-1)^{|w|(|u|+|v|)}[v,u]\c\a(w)\nonumber\\&&+[\a(w),v\c
u]-(-1)^{|u||v|}\big([w,u]\c\a(v)+(-1)^{|v|(|u|+|w|)}\a(v)\c[w,u]+[w\c
u,\a(v)]\big),\label{ee5}\\
0&=&\big(2(w\c v)\c\a(u)+(-1)^{|u|(|v|+|w|)}\a(u)\c(w\c
v)\big)-(-1)^{|u||v|}\a(w)\c(u\c v)\nonumber\\&&-(-1)^{|u||v|}\big((w\c
u)\c\a(v)+(-1)^{|u||w|}(u\c
w)\c\a(v)\big),\label{ee6}\\
0&=&(w\c v)\c\a(u)+(-1)^{|v||w|}(v\c w)\c\a(u)+\a(w)\c(v\c u)\nonumber\\&&-(-1)^{|u||v|}\big(2(w\c
u)\c\a(v)+(-1)^{|v|(|u|+|w|)}\a(v)\c(w\c u)\big),\label{ee7}\\
0&=&(w\c v)\c\a(u)+(-1)^{|u|(|v|+|w|)}\a(u)\c(w\c v)\nonumber\\&&-(-1)^{|u||v|}\big(\a(w)\c(u\c v)+(-1)^{|w|(|u|+|v|)}(u\c v)\c\a(w)\big),\label{ee8}\\
0&=&(-1)^{|u||v|}\big((w\c u)\c\a(v)+(-1)^{|v|(|u|+|w|)}\a(v)\c(w\c u)\big)-\a(w)\c(v\c u)\nonumber\\&&-(-1)^{|w|(|u|+|v|)}(v\c u)\c\a(w),\label{ee9}\\
0&=&(w\c v)\c\a(u)+(-1)^{|u|(|v|+|w|)}\a(u)\c(w\c v)+(-1)^{|v||w|}(v\c w)\c\a(u)\nonumber\\&&+(-1)^{|u|(|v|+|w|)+|v||w|}\a(u)\c(v\c w)-(-1)^{|u||v|}(w\c
u)\c\a(v)\nonumber\\&&-(-1)^{|v||w|}\a(v)\c(w\c u)-(-1)^{|u|(|w|+|v|)}(u\c
w)\c\a(v)-(-1)^{(|v|+|u|)|w|}\a(v)\c(u\c w)\nonumber\\&&-\big((-1)^{|u||v|}\a(w)\c (u\c v)+(-1)^{|w|(|u|+|v|)+|u||v|}(u\c
v)\c\a(w)\nonumber\\&&-\a(w)\c(v\c u)-(-1)^{|w|(|u|+|v|)}(v\c u)\c\a(w)\big).\label{ee10}
\end{eqnarray}
Relations \eqref{ee1}--\eqref{ee3} hold due to \eqref{super-JJ}, \eqref{HGB} and \eqref{Hom-Nov2++}, respectively. Since $u,v,w$ are arbitrary, \eqref{ee4}, \eqref{ee6} and \eqref{ee8} are equivalent to \eqref{ee5}, \eqref{ee7} and \eqref{ee9}, respectively. Rewritte \eqref{ee4} as
\begin{eqnarray*}
0&=&\big([w\c v,\a(u)]+[w,v]\c\a(u)-\a(w)\c[v,u]-(-1)^{|u||v|}[w\c
u,\a(v)]-(-1)^{|u||v|}[w,u]\c\a(v)\big)\\&&
+(-1)^{|u||v|+|u||w|+|v||w|}\big([u\c v,\a(w)]-(-1)^{|v||w|}[u\c
w,\a(v)]-(-1)^{|v||w|}[u,w]\c\a(v)\\&& +[u,v]\c\a(w)-\a(u)\c[v,w]\big),
\end{eqnarray*}
which is implied by \eqref{HGB}. By \eqref{Hom-Nov2++},  \eqref{ee6} amounts
to
\begin{eqnarray*}
(w\c v)\c\a(u)-\a(w)\c(v\c u)=(-1)^{|v||w|}\big((v\c w)\a(u)-\a(v)\c(w\c u)\big),
\end{eqnarray*}
which holds due to \eqref{Hom-Nov1++}. Relation \eqref{ee8} can be similarly obtained. Finally, rewrite \eqref{ee10} as
\begin{eqnarray*}
0&=&(-1)^{|v||w|}\big((v\c w)\c\a(u)-\a(v)\c(w\c u)-(-1)^{|v||w|}((w\c v)\c\a(u)-\a(w)\c(v\c
u))\big)\\&&+(-1)^{|u||v|}\big((w\c u)\c\a(v)-\a(w)\c(u\c v)-(-1)^{|u||w|}((u\c
w)\c\a(v)-\a(u)\c(w\c v))\big)\\&&  -(-1)^{|u||v|+|u||w|+|v||w|}\big((u\c
v)\c\a(w)-\a(u)\c(v\c w)\nonumber\\&& -(-1)^{|u||v|}(v\c u)\c\a(w)+(-1)^{|u||v|}\a(v)\c(u\c w)\big),
\end{eqnarray*}
which follows from \eqref{Hom-Nov1++}. Thus the Hom-Jacobi identity \eqref{HL3} is satisfied and $(\mathcal{R},[\d_\l\d],\a)$ is a Hom-Lie
conformal superalgebra, which is obviously quadratic by \eqref{**}.

Conversely, if $(\mathcal{R},[\d_\l\d],\a)$ is a Hom-Lie conformal superalgebra, then the Hom-Jacobi identity \eqref{HL3} gives \eqref{ee1}--\eqref{ee10}, which imply \eqref{super-ss}, \eqref{super-JJ}, \eqref{Hom-Nov1++}, \eqref{Hom-Nov2++} and \eqref{HGB} by the discussions above. Therefore, $(\mathcal {A},[\d,\d],\circ,\a)$ is a super Hom-GD bialgebra.
 \QED\vskip8pt

We call $(\mathcal{R},[\d_\l\d],\a)$ from Theorem \ref{main} is the quadratic Hom-Lie
conformal superalgebra associated to the super Hom-GD bialgebra  $(\mathcal {A},[\d,\d],\circ,\a)$.

 \begin{exam}\rm Let $(\A,[\cdot,\cdot],\circ,\alpha)$ be the three-dimensional super Hom-GD bialgebra constructed in Example \ref{exam2}. By Theorem \ref{main}, associated to $(\mathcal {A},[\d,\d],\circ,\a)$, there is a quadratic Hom-Lie
conformal superalgebra $(\mathcal{R},[\d_\l\d],\a)$, in which $\mathcal{R}=\mathbb{C}[\p]\A$ is a free $\C[\p]$-module over $\A$, $\a$ is defined by \eqref{*}, and the $\lambda$-bracket is determined by
\begin{eqnarray*}
&&[x_{1\,\l} x_1]=(\p+2\l)x_2,\ \ \ [x_{1\,\l}y]=(\l-1)y,\ \ \ \,[x_{1\,\l} x_2]=0,\\
&&[x_{2\,\l} x_2]=(\p+2\l)x_1,\ \ \ [x_{2\,\l}y]=(\l-1)y,\ \ \ \,[x_{2\,\l} x_1]=0,\\
&&[y_\l x_1]=(\p+\l+1) y,\ \ \, [y_\l x_2]=(\p+\l+1) y,\ \ [y_\l y]=0.
\end{eqnarray*}
\end{exam}

\vskip8pt

\cl{\bf\S5. \ Central extensions
}\setcounter{section}{5}\setcounter{equation}{0}\setcounter{theo}{0}
\vskip8pt
In this section, we discuss central extensions and in
particular one-dimensional central extensions of Hom-Lie conformal superalgebras. We characterize one-dimensional central extensions of quadratic Hom-Lie conformal superalgebras by using bilinear forms of the corresponding super Hom-Gel'fand-Dorfman bialgebras, generalizing a result due to \cite{H}.

Let us start with extending the notion of extension of Lie conformal superalgebras to Hom-Lie conformal superalgebras.

\begin{defi}{\rm
An {\it extension} of a Hom-Lie conformal superalgebra $(\mathcal {A},[\cdot_\l \cdot],\a)$ by an
abelian Hom-Lie conformal superalgebra $(\mathfrak a, [\cdot_\l \cdot],\a_{\mathfrak a}) $
is a commutative diagram with exact rows
\begin{equation*}
\begin{CD}
0@>>> {\mathfrak a} @>{\rm \iota}>{\rm }> {\hat \A} @>{\rm pr}>> \A @>>> 0\\
@. @V{\a_{\mathfrak a}}VV @V{\hat \a}VV @V{\a}VV\\
0@>>> {\mathfrak a} @>{\rm \iota}>{\rm }> {\hat \A} @>{\rm pr}>>
\A @>>> 0
\end{CD}
\end{equation*}
where $\rm pr$ and $\iota$ are the natural projection and
inclusion, respectively, i.e.,
\begin{eqnarray*}
\!\!\!\!\!\!\!\!\!\!\!\!&&
{\rm pr}:\hat \A\rightarrow \A,\qquad {\rm pr}(x+a)=x, \ \forall \ x\in\A, \ a\in\mathfrak{a};\\
\!\!\!\!\!\!\!\!\!\!\!\!&& \iota:\mathfrak{a}\rightarrow \hat \A,\ \, \qquad
\iota(a)=0+a, \ \forall \ a\in\mathfrak{a}.
\end{eqnarray*}
In this case, $\hat{ \mathcal{A}}$ is also called an extension of
$\mathcal{A}$ by $\mathfrak a$. The extension is said to be {\it
central} if $\p(\mathfrak a)=0$, and
\begin{eqnarray*}
{\mathfrak a}\subseteq Z(\hat{\mathcal{A} })=\{x\in \hat{\mathcal{A}
}\,\big|\, \widehat{[x_\lambda y ]}=0| \, \forall \ y\in \hat{\mathcal{A}}\}.
\end{eqnarray*}}
\end{defi}

In the following we focus on the central extension $\hat
{\mathcal{A}}$ of $\mathcal{A}$ by a one-dimensional center $\mathbb{C}
{\mathfrak c}$. This means that $\hat {\mathcal{A}}\cong
\mathcal{A}\oplus \mathbb{C} {\bf {\mathfrak c}}$, such that $\p(\mathfrak c)=0$ with $|c|=0$ and
the $\l$-bracket on $\hat{\mathcal{A}}$ is of the form
\begin{eqnarray*}
\widehat{[a_\lambda b]}=[a_\lambda
b]+f_\lambda(a,b){\mathfrak c}, \ \forall \ a,b\in\A,
\end{eqnarray*}
where $f_\lambda:\mathcal{A}\times \mathcal{A}\rightarrow \mathbb{C}[\lambda]$ is a
bilinear map. It follows from Definition \ref{HLCSA}
that $f_\lambda$ satisfies the following three properties for all homogenous elements $u, v, w\in \A$:
\begin{eqnarray}
f_\lambda(u,v)&=&-(-1)^{|u||v|}f_{-\lambda-\p}(v,u),\label{a18}\\
f_\lambda(\p u,v)&=&-\lambda f_\lambda(u,v)=-f_\lambda(u,\p v),\label{a19}\\
f_{\lambda+\mu}([u_\lambda v],\a(w))&=&f_\lambda(\a(u),[v_\mu w])-(-1)^{|u||v|}f_\mu(\a(v),[u_\lambda w])\label{a20}.
\end{eqnarray}
 The map $f_\lambda$ satisfying
\eqref{a18}--\eqref{a20} is called a {\it $2$-cocycle}.

Let $(\mathcal{R},[\d_\l\d],\a)$ be a quadratic Hom-Lie
conformal superalgebra associated to a super Hom-GD bialgebra  $(\mathcal {A},[\d,\d],\circ,\a)$. Consider the one-dimensional central extension $\hat {\R}=\R\oplus \mathbb{C} {\mathfrak c}$. According to Theorem \ref{main}, the
$\lambda$-bracket on $\hat {\R}$ is determined by
\begin{eqnarray}\label{7**}
\widehat{[u_{\lambda}v]}=[v,u]+(\partial+\lambda)(v\c
u)+(-1)^{|u||v|}\lambda(u\c v)+f_\l(u,v){\mathfrak c}, \ \forall \ u, v\in\A,
\end{eqnarray}
where $f_\lambda$ is a $2$-cocycle.

The following result is a Hom-analogue of \cite[Theorem 3.1]{H}.

\begin{theo} With notations above. Assume that $f_\l(u,v)=\mbox{$\sum^n_ {i=0}$}\l^i f_i(u,v)$ with $f_n(u,v)\neq 0$. For all homogenous elements $u, v, w\in\A$, we have
\begin{itemize}
\item[\rm(1)] If $n>3$, $f_n(u\circ v,\a(w))=0.$
\item[\rm(2)] If $n\leq 3$, then
\begin{eqnarray}
f_i(u,v)=(-1)^{|u||v|+i+1}f_i(v,u),\ \mbox{for} \ i= 0,\cdots,3;\label{ex4}
\end{eqnarray}
and
\begin{eqnarray}
&&(-1)^{|u||v|}f_3(u\circ v,\a(w))=f_3(\a(u),w\circ v)=f_3(v\circ u,\a(w)), \label{ex5}\\
&&f_3([v,u],\a(w))+(-1)^{|u||v|}f_2(u\circ v,\a(w))=f_3(\a(u),[w,v])+f_2(\a(u),w\circ v),\label{ex6}\\
&&3f_3([v,u],\a(w))-2f_2(v\circ u,\a(w))+(-1)^{|u||v|}f_2(u\circ v,\a(w))\nonumber\\
&&\ \ \ \ \ =-(-1)^{|u||v|}f_2(\a(v),w\circ u+(-1)^{|u||w|}u\circ w),\label{ex7}\\
 &&f_2([v,u],\a(w))+(-1)^{|u||v|}f_1(u\circ v,\a(w))=f_2(\a(u),[w,v])+f_1(\a(u),w\circ v),\label{ex8}\\
 &&f_1(\a(u),w\circ v+(-1)^{|v||w|}v\circ w)-(-1)^{|u||v|}f_1(\a(v),w\circ u+(-1)^{|u||w|}u\circ w)\nonumber\\
 &&\ \ \ \ \ =2f_2([v,u],\a(w))-f_1(v\circ u,\a(w))+(-1)^{|u||v|}f_1(u\circ v,\a(w)),\label{ex9}\\
 &&f_1(\a(u),[w,v])+f_0(\a(u),w\circ v)
 -(-1)^{|u||v|}f_0(\a(v),w\circ u+(-1)^{|u||w|}u\circ w)
 \nonumber\\&&\ \ \ \ \  =f_1([v,u],\a(w))+(-1)^{|u||v|}f_0(u\circ v,\a(w)),\label{ex10}\\
 &&f_0(\a(u),[w,v])-(-1)^{|u||v|}f_0(\a(v),[w,u])=f_0([v,u]\a(w))\label{ex11}.
\end{eqnarray}
\end{itemize}
\end{theo}
\noindent{\it Proof.~} By $\p(\mathfrak c)=0$, relation \eqref{a18} is written as $f_\lambda(u,v)=-(-1)^{|u||v|}f_{-\lambda}(v,u)$. This amounts to
$\mbox{$\sum^n_ {i=0}$}\l^i f_i(u,v)=-(-1)^{|u||v|}\mbox{$\sum^n_ {i=0}$}(-\l)^i f_i(v,u),$ in which
equating the coefficients of $\l^i$ gives
\begin{eqnarray}\label{k2}
f_i(u,v)=(-1)^{|u||v|+i+1}f_i(v,u), \ \mbox{for} \ \ i=0,1,\cdots,n.
\end{eqnarray}
By \eqref{a18}, \eqref{a19} and \eqref{7**}, relation \eqref{a20} is rewritten as
\begin{eqnarray}\label{k1}
&&f_{\lambda+\mu}([v,u],\a(w))-\mu f_{\lambda+\mu}(v\circ u,\a(w))+(-1)^{|u||v|}\lambda f_{\lambda+\mu}(u\circ v,\a(w))\nonumber\\
&=&f_\lambda(\a(u),[w,v])+\lambda f_\lambda(\a(u),w\circ v)+\mu f_\lambda(\a(u),w\circ v+(-1)^{|v||w|}v\circ w)\nonumber\\
&&-(-1)^{|u||v|}\big(f_\mu(\a(v),[w,u])+(\mu+\lambda) f_\mu(\a(v),w\circ u)+(-1)^{|u||w|}\lambda f_\mu(\a(v),u\circ w)\big).
\end{eqnarray}
Assume that $n>3$. Plugging $f_\l(u,v)=\mbox{$\sum^n_ {i=0}$}\l^if_i(u,v)$ into \eqref{k1}
and comparing the coefficients of $\l^2\mu^{n-1}$ and $\l^{n-1}\mu^2$, respectively, we obtain
\begin{eqnarray*}
&&-\mbox{${n\choose 2}$}f_n(v\circ u,\a(w))+(-1)^{|u||v|}nf_n(u\circ v,\a(w))=0,\\
&& -nf_n(v\circ u,\a(w))+(-1)^{|u||v|}\mbox{${n\choose 2}$}f_n(u\circ v,\a(w))=0.
\end{eqnarray*}
It follows that $f_n(u\circ v,\a(w))=0$ when $n>3.$ This proves the first assertion.

To prove the second assertion, we assume $n\leq 3$. Substituting $f_\l(u,v)=\mbox{$\sum^3_ {i=0}$}\l^i f_i(u,v)$ into \eqref{k1} and comparing the coefficients of $\l^4, \l^3\mu, \l^2\mu^2, \l\mu^3$ and $\mu^4$, respectively, we obtain
 \begin{eqnarray}
&& f_3(\a(u),w\circ v)=(-1)^{|u||v|}f_3(u\circ v,\a(w)),\label{f3-1}\\
&&f_3(\a(u),w\circ v+(-1)^{|v||w|}v\circ w)=-f_3(v\circ u,\a(w))+(-1)^{|u||v|}3f_3(u\circ v,\a(w)),\label{f3-2}\\
&&f_3(v\circ u,\a(w))=(-1)^{|u||v|}f_3(u\circ v,\a(w)),\label{f3-3}\\
&& f_3(-3v\circ u+(-1)^{|u||v|}u\circ v,\a(w))=-(-1)^{|u||v|}f_3(\a(v),w\circ u+(-1)^{|u||w|}u\circ w),\label{f3-4}\\
&&f_3(v\circ u,\a(w))=(-1)^{|u||v|}f_3(\a(v),w\circ u).\label{f3-5}
\end{eqnarray}
Because $u, v, w$ are arbitrary, \eqref{f3-1} is equivalent to \eqref{f3-5}, whereas \eqref{f3-2} amounts to \eqref{f3-4}.
Then \eqref{f3-1} and \eqref{f3-3} gives \eqref{ex5}. By \eqref{f3-1} and \eqref{f3-3} again, \eqref{f3-2} is reduced to
\begin{eqnarray}\label{f3-6}
f_3(\a(u),w\circ v)=(-1)^{|v||w|}f_3(\a(u),v\circ w).
\end{eqnarray}
By \eqref{ex4}, \eqref{f3-6} amounts to $f_3(w\circ v,\a(u))=(-1)^{|v||w|}f_3(v\circ w,\a(u))$, which is equivalent to
\eqref{f3-3}. Thus \eqref{ex5} follows from \eqref{f3-1}--\eqref{f3-5} and \eqref{ex4} when $n=3$.

Equating the coefficients of $\l^3, \l^2 \mu, \l\mu^2$ and $\mu^3$ respectively, we have
\begin{eqnarray}
 &&f_3([v,u],\a(w))+(-1)^{|u||v|}f_2(u\circ v,\a(w))=f_3(\a(u),[w,v])+f_2(\a(u),w\circ v),\label{f2-1}\\
 &&3f_3([v,u],\a(w))-f_2(v\circ u,\a(w))+(-1)^{|u||v|}2f_2(u\circ v,\a(w))\nonumber\\
 &&\ \ \ \ \ \ \ \ \ \ \ \ \ \ \ \ \ \ \ \ \ \,
 =f_2(\a(u),w\circ v+(-1)^{|v||w|}v\circ w),\label{f2-2}\\
 &&3f_3([v,u],\a(w))-2f_2(v\circ u,\a(w))+(-1)^{|u||v|}f_2(u\circ v,\a(w))\nonumber\\
 &&\ \ \ \ \ \ \ \ \ \ \ \ \ \ \ \ \ \ \ \ \ \,  =-(-1)^{|u||v|}f_2(\a(v),w\circ u+(-1)^{|u||w|}u\circ w), \label{f2-3}\\
 &&f_3([v,u],\a(w))-f_2(v\circ u,\a(w))=-(-1)^{|u||v|}(f_3(\a(v),[w,u])+f_2(\a(v),w\circ u)).\label{f2-4}
 \end{eqnarray}
Letting $u$ and $v$ change places in \eqref{f2-4} gives
\begin{eqnarray*}
f_3([u,v],\a(w))-f_2(u\circ v,\a(w))=-(-1)^{|u||v|}(f_3(\a(u),[v,u])+f_2(\a(u),w\circ v)).
\end{eqnarray*}
It is equivalent to \eqref{f2-1}, which is exactly \eqref{ex6}. By similar discussions, \eqref{f2-2} is equivalent to \eqref{f2-3}. Thus we get \eqref{ex7}.

Finally, extracting the coefficients of $\l^2, \l \mu, \mu^2, \l,\mu$ and $\l^0\mu^0$, respectively, gives
 \begin{eqnarray}
 &&f_2([v,u],\a(w))+(-1)^{|u||v|}f_1(u\circ v,\a(w))=f_2(\a(u),[w,v])+f_1(\a(u),w\circ v),\label{f1-1}\\
 &&f_1(\a(u),w\circ v+(-1)^{|v||w|}v\circ w)
  -(-1)^{|u||v|}f_1(\a(v),w\circ u+(-1)^{|u||w|}u\circ w)\nonumber\\
  &&\ \ \ \ \ \ \ \ \ \ \ \ \ \ \ \ \ \ \ \ =2f_2([v,u],\a(w))-f_1(v\circ u,\a(w))+(-1)^{|u||v|}f_1(u\circ v,\a(w)),\label{f1-2}\\
 &&f_2([v,u],\a(w))-f_1(v\circ u,\a(w))=-(-1)^{|u||v|}(f_2(\a(v),[w,u])-f_1(\a(v),w\circ u)),\label{f1-3}\\
 &&f_1([v,u],\a(w))+(-1)^{|u||v|}f_0(u\circ v,\a(w))=f_1(\a(u),[w,v])+f_0(\a(u),w\circ v)\nonumber\\
 &&\ \ \ \ \ \ \ \ \ \ \ \ \ \ \ \ \ \ \ \ -(-1)^{|u||v|}f_0(\a(v),w\circ u+(-1)^{|u||w|}u\circ w),\label{f1-4}\\
 &&f_1([v,u],\a(w))-f_0(v\circ u,\a(w))=f_0(\a(u),w\circ v+(-1)^{|v||w|}v\circ w)\nonumber\\
 &&\ \ \ \ \ \ \ \ \ \ \ \ \ \ \ \ \ \ \ \ -(-1)^{|u||v|}f_1(\a(v),[w,u])-(-1)^{|u||v|}f_0(\a(v),w\circ u),\label{f1-5}\\
 &&f_0(\a(u),[w,v])-(-1)^{|u||v|}f_0(\a(v),[w,u])=f_0([v,u],\a(w))\label{f0}.
 \end{eqnarray}
We see that \eqref{f1-2} is exactly \eqref{ex9}, whereas \eqref{f0} is exactly \eqref{ex11}. As discussed in \eqref{f2-4},
\eqref{f1-1} and \eqref{f1-4} are equivalent to \eqref{f1-3} and \eqref{f1-5}, respectively. Thus we obtain \eqref{ex8} and
\eqref{ex10}. This completes the proof.\QED

 \vskip10pt

\small\noindent{\bf Acknowledgements~}~{\footnotesize This work was supported by National Natural Science
Foundation grants of China (11301109), the Research Fund for the Doctoral Program of Higher Education
(20132302120042), and the Fundamental Research
Funds for the Central Universities (HIT. NSRIF. 201462).}\\

\end{document}